\documentclass[11pt]{article}

\usepackage[utf8]{inputenc}
\usepackage{amssymb,epsfig}
\usepackage{amsmath}
\usepackage{comment}
\usepackage[dvips]{color}
\usepackage[T1]{fontenc}

\newtheorem{defi}{Definition}
\newtheorem{prop}[defi]{Proposition}
\newtheorem{theo}[defi]{Theorem}
\newtheorem{conj}[defi]{Conjecture}
\newtheorem{lemm}[defi]{Lemma}
\newtheorem{coro}[defi]{Corollary}
\newtheorem{rema}[defi]{Remark}
\newtheorem{exem}[defi]{Example}
\newtheorem{exems}[defi]{Examples}

\newcommand{\bdefi}{\begin{defi}}
\newcommand{\edefi}{\end{defi}}
\newcommand{\bprop}{\begin{prop}}
\newcommand{\eprop}{\end{prop}}
\newcommand{\btheo}{\begin{theo}}
\newcommand{\etheo}{\end{theo}}
\newcommand{\blemm}{\begin{lemm}}
\newcommand{\brema}{\begin{rema}}
\newcommand{\erema}{\end{rema}}
\newcommand{\bexer}{\begin{exem}}
\newcommand{\eexer}{\end{exem}}
\newcommand{\bexems}{\begin{exems}}
\newcommand{\eexems}{\end{exems}}
\newcommand{\bconj}{\begin{conj}}
\newcommand{\econj}{\end{conj}}
\newcommand{\elemm}{\end{lemm}}
\newcommand{\bcoro}{\begin{coro}}
\newcommand{\ecoro}{\end{coro}}
\newcommand{\dem}{\noindent{\bf Proof. }}
\newcommand{\rem}{\noindent{\bf Remark. }}

\usepackage{mathrsfs}
\renewcommand\mathcal{\mathscr}

\newcommand{\M}{{\cal M}}
\newcommand{\N}{{\cal N}}
\newcommand{\G}{{\cal G}}
\newcommand{\D}{{\cal D}}
\newcommand{\E}{{\cal E}}
\newcommand{\F}{{\cal F}}

\renewcommand{\H}{{\cal H}}
\newcommand{\OOO}{{\cal O}}
\newcommand{\C}{{\cal C}}
\newcommand{\I}{{\cal I}}
\newcommand{\Q}{{\cal Q}}
\renewcommand{\P}{{\cal P}}

\newcommand{\maths}[1]{{\mathbb #1}}  

\renewcommand{\AA}{\maths{A}}

\newcommand{\RR}{\maths{R}}
\newcommand{\NN}{\maths{N}}
\newcommand{\CC}{\maths{C}}
\newcommand{\QQ}{\maths{Q}}

\newcommand{\SSS}{\maths{S}}

\newcommand{\HH}{\maths{H}}
\newcommand{\FF}{\maths{F}}

\newcommand{\ZZ}{\maths{Z}}
\newcommand{\PP}{\maths{P}}

\newcommand{\aaa}{{\mathfrak a}}

\newcommand{\uuu}{{\mathfrak u}}
\newcommand{\vvv}{{\mathfrak v}}

\newcommand{\mmm}{{\mathfrak m}}
\renewcommand{\uuu}{{\mathfrak u}}

\newcommand{\ra}{\rightarrow}
\newcommand{\bs}{\backslash}

\newcommand{\ov}[1]{{\overline #1}} 
\newcommand{\wt}[1]{{\widetilde{#1}}}
\newcommand{\wh}[1]{{\widehat{#1}}}

\newcommand{\ga}{\gamma}
\newcommand{\Ga}{\Gamma}

\newcommand{\cqfd}{\hfill$\Box$}

\newcommand{\card}{{\operatorname{Card}}}
\renewcommand{\Re}{{\operatorname{Re}}}

\newcommand{\Vol}{\operatorname{Vol}}

\newcommand{\covol}{\operatorname{Covol}}
\newcommand{\id}{\operatorname{id}}
\newcommand{\PSL}{\operatorname{PSL}}
\newcommand{\SL}{\operatorname{SL}}
\newcommand{\Res}{\operatorname{Res}}
\newcommand{\dvol}{\;d\operatorname{vol}}
\newcommand{\Det}{\operatorname{Det}}

\newcommand{\hcr}{{\HH}^5_\RR}
\newcommand{\hnr}{{\HH}^n_\RR}
\newcommand{\phcr}{\partial_\infty\hcr}

\newcommand{\SLO}{\operatorname{SL}_{2}(\OOO)}
\newcommand{\GLO}{\operatorname{GL}_{2}(\OOO)}
\newcommand{\PSLO}{\operatorname{PSL}_{2}(\OOO)}
\newcommand{\SLH}{\operatorname{SL}_{2}(\HH)}
\newcommand{\PSLH}{\operatorname{PSL}_{2}(\HH)}
\newcommand{\autom}{\operatorname{SU}_f(\OOO)}
\newcommand{\Pautom}{\operatorname{PSU}_f(\OOO)}

\newcommand{\SpO}{\operatorname{Sp}_1(\OOO)}

\newcommand{\tr}{\operatorname{tr}}
\newcommand{\n}{\operatorname{n}}

\newcommand\Rep{\mathcal R}

\newcounter{fig}

\title{On the arithmetic and geometry \\of binary Hamiltonian forms}
\author{Jouni Parkkonen \and Fr\'ed\'eric Paulin} 
\date{With an appendix by Vincent Emery}

\begin{document}
\bibliographystyle{../alphanum}

\maketitle

\begin{abstract} 
  Given an indefinite binary quaternionic Hermitian form $f$ with
  coefficients in a maximal order of a definite quaternion algebra
  over $\QQ$, we give a precise asymptotic equivalent to the number of
  nonequivalent representations, satisfying some congruence
  properties, of the rational integers with absolute value at most $s$
  by $f$, as $s$ tends to $+\infty$. We compute the volumes of
  hyperbolic $5$-manifolds constructed by quaternions using Eisenstein
  series. In the Appendix, V.~Emery computes these volumes using
  Prasad's general formula. We use hyperbolic geometry in dimension
  $5$ to describe the reduction theory of both definite and indefinite
  binary quaternionic Hermitian forms.  
\end{abstract}

\section{Introduction}
\label{sect:intro}

Following H.~Weyl \cite{Weyl40}, we will call {\it Hamiltonian
  form} a Hermitian form over Hamilton's real quaternion algebra with
anti-involution the conjugation.

Since Gauss, the reduction theory of the integral binary quadratic
forms, and the problem of representation of integers by them, is quite
completely understood. For binary Hermitian forms, these subjects have
been well studied, starting with Hermite, Bianchi and especially
Humbert, and much developed by Elstrodt, Grunewald and Mennicke (see
for instance \cite{ElsGruMen98}). In the recent paper
\cite{ParPau11CJM}, we gave a precise asymptotic on the number of
nonequivalent proper representations of rational integers with
absolute value at most $s$ by a given integral indefinite Hermitian
form.  Besides the general results on quadratic forms (see for
instance \cite{Weyl40,Cassels08}) and some special works (see for
instance \cite{Pronin67}, \cite{HasIbu80}), not much seemed to be
precisely known on these questions for binary Hamiltonian forms.

In this paper, we use hyperbolic geometry in dimension $5$ to study
the asymptotic of the counting of representations of rational integers
by binary Hamiltonian forms and to give a geometric description of the
reduction theory of such forms. General formulas are known (by
Siegel's mass formula, see for instance \cite{EskRudSar91}), but it
does not seem to be easy (or even doable) to deduce our asymptotic
formulas from them.  There are numerous results on the counting of
integer points with bounded norm on quadrics (or homogeneous
varieties), see for instance the works of Duke, Eskin, McMullen, Oh,
Rudnick, Sarnak and others. In this paper, we count appropriate orbits
of integer points on which a fixed integral binary Hamiltonian form is
constant, analogously to \cite{ParPau11CJM}.

\medskip Let $\HH$ be Hamilton's quaternion algebra over $\RR$, with
$x\mapsto \overline{x}$ its conjugation, $\n: x\mapsto x\overline{x}$
its reduced norm and $\tr: x\mapsto x+\overline{x}$ its reduced
trace. Let $A$ be a quaternion algebra over $\QQ$, which is definite
($A\otimes_\QQ\RR=\HH$), with reduced discriminant $D_A$ and class
number $h_A$. Let $\OOO$ be a maximal order in $A$, and let $\mmm$ be
a (nonzero) left fractional ideal of $\OOO$, with reduced norm
$\n(\mmm)$ (see Section \ref{sect:defquadalg} for definitions).

Let $f:\HH\times\HH\ra \RR$ be a binary Hamiltonian form,
with
\begin{equation}\label{eq:form}
f(u,v)=a\,\n(u)+ \tr(\ov u\, b \,v) +c\,\n(v)\;,
\end{equation}
which is integral over $\OOO$ (its coefficients satisfy $a,c\in\ZZ$
and $b\in\OOO$) and indefinite (its discriminant $\Delta(f)=\n(b)-ac$
is positive), see Section \ref{sect:indefbinhamform}.  We denote by
$\SLO$ the group of invertible $2\times 2$ matrices with coefficients
in $\OOO$ (see Section \ref{sect:hambiagroup}).  The group $\autom$ of
automorphs of $f$ consists of those elements $g\in\SLO$ for which
$f\circ g=f$. Given an arithmetic group $\Ga$, such as $\SLO$ and
$\autom$, we will denote by $\covol(\Ga)$ the volume of the quotient
by $\Ga$ of its associated symmetric space (assumed to be of
noncompact type and normalized to have $-1$ as
the minimum of its sectional curvature).

For every $s>0$, we consider the integer
$$
\psi_{f,\mmm}(s)=\card\;\;_{\mbox{$\autom$}}\bs
\big\{(u,v)\in\mmm\times\mmm\;:\;\n(\mmm)^{-1}|f(u,v)|\leq s,
\;\;\;\OOO u+\OOO v=\mmm\big\}\;,
$$
which is the number of nonequivalent $\mmm$-primitive representations
by $f$ of rational integers with absolute value at most $s$.  The
finiteness of $\psi_{f,\mmm}(s)$ follows from general results on
orbits of algebraic groups defined over number fields
\cite[Lem.~5.3]{BorelHarishChandra62}.

\btheo\label{theo:mainintro} As $s$ tends to $+\infty$, we have the
equivalence, with $p$ ranging over positive rational primes,
$$
\psi_{f,\mmm}(s)\sim \frac{45\;D_A\;\covol(\autom)}
{2\pi^2\;\zeta(3)\;\Delta(f)^2\;\prod_{p|D_A}(p^3-1)}\;\;s^4\;.
$$
\etheo 

This result follows from the more general Theorem
\ref{theo:mainversionG}, which allows in particular to count
representations satisfying given congruence properties (see the end of
Section \ref{sect:rephamil}).

Here is an example of our applications, concerning the asymptotic of
the very useful real scalar product $(u,v)\mapsto \tr
(\overline{u}\,v)$ on $\HH$. See Section \ref{sect:rephamil} for the
proof and for further applications.  Let
$$
\SpO=\Big\{g\in\SLO\;:\;
^t\overline{g}\,\Big(\begin{array}{cc}0& 1\\1&
  0\end{array}\Big)\,g=\Big(\begin{array}{cc}0& 1\\1&
  0\end{array}\Big)\Big\}\;.
$$

\bcoro\label{coro:introspunooo} 
As $s$ tends to $+\infty$, we have the equivalence\\

\noindent ${\displaystyle
\card\;\;_{\mbox{$\SpO$}}\bs
\big\{(u,v)\in \OOO\times \OOO\;:\;
|\tr(\ov u\,v)|\leq s,\;\OOO u+ \OOO v=\OOO\;\big\}}$

\hfill $\displaystyle
\sim \; \frac{D_A}
{48\,\zeta(3)}\prod_{p|D_A}\frac{p^2+1}{p^2+p+1}\;\;s^4\;.
$
\ecoro

In order to prove Theorem \ref{theo:mainintro}, applying a counting
result of \cite{ParPau} following from dynamical properties of the
geodesic flow of real hyperbolic manifolds, we first prove that
$$
\psi_{f,\mmm}(s)\sim \frac{D_A\prod_{p|D_A}(p-1)\,\covol(\autom)}
{512\,\pi^2\,\Delta(f)^2\,\covol(\SLO)}\;s^4\;.
$$ 
The covolumes of the arithmetic groups $\SLO$ and $\autom$ may be
computed using Prasad's very general formula in \cite{Prasad89} (see
\cite{Emery09} for an excellent exposition). Following the approach of
Rankin-Selberg \cite{Rankin39,Selberg40}, Langlands
\cite{Langlands66b}, Sarnak \cite{Sarnak83} and others, we compute
$\covol(\SLO)$ in the main body of this paper (see Section
\ref{sect:Eisenstein}) using Eisenstein series, whose analytic
properties in the quaternion setting have been studied by
\cite{KraOse90}. We initially proved the case $h_A=1$ of the following
result, V.~Emery proved the general case using Prasad's formula (see
the Appendix), and we afterwards managed to push the Eisenstein series
approach to get the general result. The two proofs are completely
different.

\btheo(Emery, see the Appendix) We have
$$
\covol(\SLO)= \frac{\zeta(3)\prod_{p|D_A}(p^3-1)(p-1)\,}{11520}\;.
$$
\etheo

In the final section, we give a geometric reduction theory of binary
Hamiltonian forms using real hyperbolic geometry. The case of binary
quadratic forms is well known, from either the arithmetic, geometric
or algorithmic viewpoint (see for instance
\cite{Cassels08,Zagier81,BucVol07}). We refer for instance to
\cite{ElsGruMen98} for the reduction theory of binary Hermitian
forms. The case of binary Hamiltonian forms has been less developed,
see for instance \cite{Pronin67,HasIbu80} for results in the positive
definite case. We construct a natural map $\Xi$ from the set
$\Q(\OOO,\Delta)$ of binary Hamiltonian forms that are integral over
$\OOO$ and have a fixed discriminant $\Delta\in\ZZ-\{0\}$ to the set
of points or totally geodesic hyperplanes of the $5$-dimensional real
hyperbolic space $\HH^5_\RR$. For $\F_\OOO$ a Ford fundamental domain
for the action of $\SLO$ on $\HH^5_\RR$, we say that
$f\in\Q(\OOO,\Delta)$ is {\it reduced} if $\Xi(f)$ meets
$\F_\OOO$. The finiteness of the number of orbits of $\SLO$ on
$\Q(\OOO,\Delta)$, which can be deduced from general results of
Borel-Harich-Chandra, then follows in an explicit way from the
equivariance property of $\Xi$ and the following result proved in
Section \ref{sect:geomredtheo}.

\btheo \label{theo:inducintro}
There are only finitely many reduced integral binary
Hamiltonian forms with a fixed nonzero discriminant.  
\etheo

Answering the remark page 257 of \cite{Cassels08} that explicit sets
of inequalities implying the reduction property were essentially only
known for quadratic forms in dimension $n\leq 7$, we give an explicit
such set in dimension $8$ at the end of Section
\ref{sect:geomredtheo}.

The knowledgeable reader may skip the background sections
\ref{sect:defquadalg} (except the new Lemma \ref{lem:chenevier}),
\ref{sect:hambiagroup} and \ref{sect:indefbinhamform} on respectively
definite quaternion algebras over $\QQ$, quaternionic homographies and
real hyperbolic geometry in dimension $5$, and binary Hamiltonian
forms, though many references are made to them in the subsequent
sections.

\medskip {\small {\it Acknowledgment. } We thank P.~Sarnak for his
  comments on the origin of volume computations using Eisenstein
  series, G.~Chenevier for the proof of Lemma \ref{lem:chenevier},
  Y.~Benoist and F.~Choucroun for discussions related to the Appendix,
  and the referee for helpful comments, in particular for Lemma
  \ref{lem:rootgroup}. The second author thanks the University of
  Jyv\"askyl\"a for the nice snow and its financial support.}

\section{Background on definite quaternion algebras 
over $\QQ$}
\label{sect:defquadalg}

A {\em quaternion algebra} over a field $F$ is a four-dimensional
central simple algebra over $F$. We refer for instance to
\cite{Vigneras80} for generalities on quaternion algebras. 

A real quaternion algebra is isomorphic either to $\M_2(\RR)$ or to
Hamilton's quaternion algebra $\HH$ over $\RR$, with basis elements
$1,i,j,k$ as a $\RR$-vector space, with unit element $1$ and
$i^2=j^2=-1$, $ij=-ji=k$. We define the {\em conjugate} of
$x=x_0+x_1i+x_2j+x_3k$ in $\HH$ by $\overline{x}=x_0-x_1i-x_2j-x_3k$,
its {\em reduced trace} by $\tr(x)=x+\overline{x}$, and its {\em
  reduced norm} by $\n(x)= x\,\overline{x}=\overline{x}\,x$. Note that
$\n(xy)=\n(x)\n(y)$, and $\n(x)\geq 0$ with equality if and only if
$x= 0$, hence $\HH$ is a division algebra, with $[\HH^\times,
\HH^\times] = \n^{-1}(1)$.  Furthermore, $\tr(\overline{x}) = \tr(x)$
and $\tr(xy)=\tr(yx)$.  For every matrix $X=(x_{i,j})_{1\leq i\leq
  p,\; 1\leq j\leq q}\in\M_{p,q}(\HH)$, we denote by $X^*=
(\overline{x_{j,i}})_{1\leq i\leq q,\;1\leq j\leq p}\in\M_{q,p}(\HH)$
its adjoint matrix, which satisfies $(XY)^*=Y^*X^*$. The matrix $X$ is
{\it Hermitian} if $X=X^*$.

Let $A$ be a quaternion algebra over $\QQ$. We say that $A$ is {\em
  definite} (or ramified over $\RR$) if the real quaternion algebra
$A\otimes_\QQ\RR$ is isomorphic to $\HH$. In this paper, whenever we
consider a definite quaternion algebra $A$ over $\QQ$, we will fix an
identification between $A\otimes_\QQ\RR$ and $\HH$, so that $A$ is a
$\QQ$-subalgebra of $\HH$.

The {\em reduced discriminant} $D_A$ of $A$ is the product of the
primes $p\in\NN$ such that the quaternion algebra $A\otimes_\QQ\QQ_p$
over $\QQ_p$ is a division algebra.  Two definite quaternion algebras
over $\QQ$ are isomorphic if and only if they have the same reduced
discriminant, which can be any product of an odd number of primes (see
\cite[page 74]{Vigneras80}).

A {\em $\ZZ$-lattice} $I$ in $A$ is a finitely generated $\ZZ$-module
generating $A$ as a $\QQ$-vector space. The intersection of finitely
many $\ZZ$-lattices of $A$ is again a $\ZZ$-lattice.  An {\it order}
in a quaternion algebra $A$ over $\QQ$ is a unitary subring $\OOO$ of
$A$ which is a $\ZZ$-lattice. In particular, $A=\QQ\OOO$. Each order
of $A$ is contained in a maximal order. The {\it type number} $t_A\geq
1$ of $A$ is the number of conjugacy (or equivalently isomorphism)
classes of maximal orders in $A$ (see for instance \cite[page
152]{Vigneras80} for a formula). For instance, $t_A=1$ if
$D_A=2,3,5,7,13$ and $t_A=2$ if $D_A=11,17$.  If $\OOO$ is a maximal
order in $A$, then the ring $\OOO$ has $2$, $4$ or $6$ invertible
elements except that $|\OOO^\times|=24$ when $D_A=2$, and
$|\OOO^\times|=12$ when $D_A=3$. When $D_A=2,3,5,7,13$, then (see
\cite[page 103]{Eichler38})
\begin{equation}\label{eq:computnombreunit}
|\OOO^\times|=\frac{24}{D_A-1}\;.
\end{equation}

\medskip
\noindent{\bf Examples.} (See \cite[page 98]{Vigneras80}.)  (1) The
$\QQ$-vector space $A=\QQ+\QQ i+\QQ j+\QQ k$ generated by $1,i,j,k$ in
$\HH$ is Hamilton's quaternion algebra over $\QQ$.  It is the unique
definite quaternion algebra over $\QQ$ (up to isomorphism) with
discriminant $D_A=2$. The {\it Hurwitz order} $\OOO=\ZZ+\ZZ i+\ZZ j+
\ZZ\frac{1+i+j+k}{2}$ is maximal, and it is unique up to conjugacy.

(2) Similarly, $A=\QQ+\QQ i+\QQ \,\sqrt{p}\,j+\QQ \,\sqrt{p}\,k$ is
the unique (up to isomorphism) definite quaternion algebra over $\QQ$
with discriminant $D_{A}=p$ for $p=3,7$, and $\OOO=\ZZ+\ZZ i+\ZZ
\frac{i+\sqrt{p}\,j}{2} + \ZZ\frac{1+\sqrt{p}\,k}{2}$ is its unique
(up to conjugacy) maximal order.

(3) Similarly, $A=\QQ+\QQ \,\sqrt{2}\,i+\QQ\, \sqrt{p}\,j+\QQ
\,\sqrt{2p}\,k$ is the unique (up to isomorphism) definite quaternion
algebra over $\QQ$ with discriminant $D_{A}=p$ for $p=5,13$, and
$\OOO=\ZZ+\ZZ \frac{1 + \sqrt{2}\,i + \sqrt{p}\,j}{2} + \ZZ
\frac{\sqrt{p}\,j}{2}+ \ZZ\frac{2+\sqrt{2}\,i +\sqrt{2p}\,k}{2}$ is
its unique (up to conjugacy) maximal order.

\bigskip Let $\OOO$ be an order in $A$. The reduced norm $\n$ and the
reduced trace $\tr$ take integral values on $\OOO$. The invertible
elements of $\OOO$ are its elements of reduced norm $1$. Since
$\overline{x}=\tr(x)-x$, any order is invariant under conjugation.

The {\em left order} $\OOO_\ell(I)$ of a $\ZZ$-lattice $I$ is $\{x\in
A\;:\; xI \subset I\}$; its {\em right order} $\OOO_r(I)$ is $\{x\in
A\;:\; Ix \subset I\}$. A {\em left fractional ideal} of $\OOO$ is a
$\ZZ$-lattice of $A$ whose left order is $\OOO$. A {\em left ideal} of
$\OOO$ is a left fractional ideal of $\OOO$ contained in $\OOO$. Right
(fractional) ideals are defined analogously.  The {\em inverse} of a
right fractional ideal $\mmm$ of $\OOO$ is $\mmm^{-1}=\{x\in
A\;:\;\mmm \,x \,\mmm\subset \mmm\}$. It is easy to check that for
every $u,v\in\OOO$, if $uv\neq 0$, then
\begin{equation}\label{eq:invsomeginter}
(u\OOO+v\OOO)^{-1}=\OOO u^{-1}\cap \OOO v^{-1}\;.
\end{equation}
If $\OOO$ is maximal, then $\mmm^{-1}$ is a left fractional ideal of
$\OOO$ and
\begin{equation}\label{eq:gauchedroite}
\OOO_r(\mmm^{-1})=\OOO_\ell(\mmm)\;.
\end{equation}
This formula follows from Lemma 4.3 (3) of \cite[page
21]{Vigneras80}, which says that $\OOO_r(\mmm^{-1})$ contains
$\OOO_\ell(\mmm)$, since the maximality of $\OOO$ implies the
maximality of $\OOO_\ell(\mmm)$, by Exercice 4.1 of \cite[page
28]{Vigneras80}.

Two  left fractional ideals $\mmm$ and $\mmm'$ of $\OOO$ are
isomorphic as left $\OOO$-modules if and only if $\mmm'=\mmm c$ for
some $c\in A^\times$. A (left) {\it ideal class} of $\OOO$ is an
equivalence class of  left fractional ideals of $\OOO$ for this
equivalence relation. We will denote by $_\OOO\!\I$ the set of ideal
classes of $\OOO$, and by $[\mmm]$ the ideal class of a  left
fractional ideal $\mmm$ of $\OOO$. The {\it class number} $h_A$ of $A$
is the number of ideal classes of a maximal order $\OOO$ of $A$. It is
finite and independent of the maximal order $\OOO$ (see for instance
\cite[page 87-88]{Vigneras80}). See for instance \cite[page
152-155]{Vigneras80} for a formula for $h_A$, and for the fact that
$h_A=1$ if and only if $D_A=2,3,5,7,13$. In particular $D_A$ is prime
if $h_A=1$.

The {\it norm} $\n(\mmm)$ of a  left (or right) ideal $\mmm$ of
$\OOO$ is the greatest common divisor of the norms of the nonzero
elements of $\mmm$. In particular, $\n(\OOO)=1$.  The {\it norm} of a
 left (or right) fractional ideal $\mmm$ of $\OOO$ is
$\frac{\n(c\mmm)}{\n(c)}$ for any $c\in \NN-\{0\}$ such that
$c\mmm\subset \OOO$.

Note that a $\ZZ$-lattice $\Lambda$ in $A$ is a $\ZZ$-lattice in the
Euclidean vector space $\HH$ (with orthonormal basis $(1,i,j,k)$), and
the volume $\Vol(\Lambda\bs\HH)$ is finite. If $\OOO$ is maximal,
we have (see for instance \cite[Lem.~5.5]{KraOse90})
\begin{equation}\label{eq:volumeHHparOOO}
\Vol(\OOO\bs\HH)=\frac{D_A}{4}\;.
\end{equation}

\medskip The classical {\it zeta function} of $A$ is 
$$
\zeta_A(s)=\sum_{\aaa}\frac{1}{\n(\aaa)^{2s}}\;,
$$ 
where the sum is over all  left ideals $\aaa$ in a maximal
order $\OOO$ of $A$.  It is independent of the choice of $\OOO$, it is
holomorphic on $\{s\in\CC\;:\;\Re\;s>1\}$ and it satisfies by a
theorem of Hey, with $\zeta$ the usual Riemann zeta function,
\begin{equation}\label{eq:valeurzetaooo}
\zeta_A(s)=\zeta(2s)\;\zeta(2s-1)\;
\prod_{p\,|\,D_A}(1-p^{1-2s})\;,
\end{equation} 
where as usual the index $p$ is prime (see \cite[page
88]{Schoeneberg39} or \cite[page 64]{Vigneras80}).
Let $\mmm$ be  a left fractional ideal of a maximal order $\OOO$ in
$A$. Define
$$
\zeta(\mmm,s)=
\n(\mmm)^{2s}\sum_{x\in \mmm-\{0\}}\;\frac{1}{\n(x)^{2s}}\;,
$$
which is also holomorphic on $\Re\;s>1$ (and depends only on the ideal
class of $\mmm$), and
$$
\zeta_{[\mmm]}(s)=\sum \frac{1}{\n(\aaa)^{2s}}
$$ 
where the sum is over all  left ideals $\aaa$ in $\OOO$ whose
ideal class is $[\mmm]$.  The relations we will use in Section
\ref{sect:Eisenstein} between these zeta functions are the following
ones, where $\Re\;s>1$. The first one is obvious, see for instance
respectively \cite[page 134]{Deuring35} and \cite[page
436]{KraOse90} for the other two:
\begin{equation}\label{eq:Deuring128}
\zeta_{A}(s)=\sum_{[\aaa]\,\in\; _\OOO\!\I}\; \zeta_{[\aaa]}(s)\;,
\end{equation} 
\begin{equation}\label{eq:Deuring134}
\sum_{[\aaa]\,\in\; _\OOO\!\I}\; \frac{1}{\big|\OOO_r(\aaa)^\times\big|} =
\frac{1}{24}\;\prod_{p|D_A} (p-1)\;,
\end{equation}
\begin{equation}\label{eq:KraOse436}
\zeta(\mmm,s)=\big|\OOO_r(\mmm)^\times\big|\;\zeta_{[\mmm]}(s)\;.
\end{equation}
Note that when the class number $h_A$ of $A$ is $1$, the formula
\eqref{eq:KraOse436} becomes
\begin{equation}\label{eq:valeurpartzetaideal}
\zeta(\OOO,s)=|\OOO^\times|\;\zeta_A(s)\;.
\end{equation} 

We end this section with the following lemma, which will be used in
the proof of Theorem \ref{theo:mainversionG}.

\blemm \label{lem:chenevier}
Let $\OOO$ be a maximal order in a definite quaternion algebra
$A$ over $\QQ$, let $z\in A-\{0\}$ and let $\Lambda=\OOO \cap z \OOO
\cap \OOO z\cap z\OOO z$. Then $\Lambda$ is a $\ZZ$-sublattice of
$\OOO$ such that
$$
[\OOO:\Lambda]\n(\OOO z^{-1}+\OOO)^4=1\;.
$$
\elemm

\dem This is a ``prime by prime'' type of proof, suggested by
G.~Chenevier. As an intersection of four $\ZZ$-lattices, $\Lambda$ is
a $\ZZ$-lattice, contained in $\OOO$. For every (positive rational)
prime $p$, let $\nu_p$ be the $p$-adic valuation on $\QQ_p$; let us
consider the quaternion algebra $A_p=A\otimes_\QQ \QQ_p$ over $\QQ_p$,
whose reduced norm is denoted by $\n_p:A_p\ra\QQ_p$; and for every
$\ZZ$-lattice $L$ of $A$, let $L_p=L\otimes_\ZZ\ZZ_p$. We embed $A$ in
$A_p$ as usual by $x\mapsto x\otimes 1$. We then have the following
properties (see for instance \cite[page 83-84]{Vigneras80}): $L_p$ is
a $\ZZ_p$-lattice of $A_p$; the map $L\mapsto L_p$ commutes with the
inclusion, the sum and the intersection; if $L,L'$ are $\ZZ$-lattices
with $L\subset L'$, then
$$
[L':L]=\prod_{p}\;[L'_p:L_p]\;;
$$
if $L$ is a left fractional ideal of
$\OOO$, then $L_p$ is a left  fractional ideal of $\OOO_p$, and
$$
\n(L)=\prod_{p}\; p^{\nu_p(\n_p(L_p))}\;.
$$
Hence in order to prove Lemma \ref{lem:chenevier}, we only have to
prove that for every prime $p$, if $z\in A_p^\times$ and $\Lambda_p=
\OOO_p\cap z\OOO_p\cap \OOO_p\,  z\cap z\OOO_p \, z$, we have
\begin{equation}\label{eq:primbyprim}
[\OOO_p:\Lambda_p]=p^{-4\,\nu_p(\n_p(\OOO_p \,z^{-1}+\OOO_p))}\;.
\end{equation} 
We distinguish two cases.

First assume that $p$ does not divide $D_A$. Then we may assume that
$A_p=\M_2(\QQ_p)$ and $\OOO_p=\M_2(\ZZ_p)$ (by the uniqueness up to
conjugacy of maximal orders). By Cartan's decomposition of
$\operatorname{GL}_2(\QQ_p)$ (see for instance \cite{BruTit72}, or
consider the action of $\operatorname{GL}_2(\QQ_p)$ on its Bruhat-Tits
tree as in \cite{Serre83}), the element
$z\in\operatorname{GL}_2(\QQ_p)$ may be written
$z=P\Big(\!\!\begin{array}{cc} p^{a} & 0 \\ 0 &
  p^{b} \end{array}\!\!\Big)Q$ with $P,Q$ in the (good) maximal
compact subgroup $\operatorname{GL}_2(\ZZ_p)$ and $a, b$ in $\ZZ$.
Since $\operatorname{GL}_2(\ZZ_p)$ preserves $\OOO_p=\M_2(\ZZ_p)$ by
left and right multiplication, preserves the indices of
$\ZZ$-lattices, and contains only elements of reduced norm (that is of
determinant) having valuation $0$, we may assume that $P=Q=\id$. We
hence have, by an easy matrix computation,
\begin{align*}
\Lambda_p&=
\Big(\!\!\begin{array}{cc} \ZZ_p\cap p^a\,\ZZ_p\cap p^{2a}\,\ZZ_p & 
\ZZ_p\cap p^a\,\ZZ_p\cap p^b\,\ZZ_p\cap p^{a+b}\,\ZZ_p 
\\ \ZZ_p\cap p^a\,\ZZ_p\cap p^b\,\ZZ_p\cap p^{a+b}\,\ZZ_p  &
\ZZ_p\cap p^b\,\ZZ_p\cap p^{2b}\,\ZZ_p \end{array}\!\!\Big)\\ &=
\Big(\!\!\begin{array}{cc} p^{2\max\{a,0\}}\,\ZZ_p & 
p^{\max\{a,0\}+\max\{b,0\}}\,\ZZ_p \\
p^{\max\{a,0\}+\max\{b,0\}}\,\ZZ_p
& p^{2\max\{b,0\}}\,\ZZ_p\end{array}\!\!\Big)\;.
\end{align*}
Similarly, we have
$$
\OOO_p\, z^{-1}+\OOO_p=
\Big(\!\!\begin{array}{cc} p^{-a}\,\ZZ_p+\ZZ_p & p^{-b}\,\ZZ_p+\ZZ_p 
\\ p^{-a}\,\ZZ_p+\ZZ_p & p^{-b}\,\ZZ_p+\ZZ_p  \end{array}\!\!\Big)
=
\M_2(\ZZ_p)\Big(\!\!\begin{array}{cc} p^{\min\{-a,0\}} & 0
\\ 0 & p^{\min\{-b,0\}} \end{array}\!\!\Big)\;.
$$
Therefore, since $\n_p(\M_2(\ZZ_p))=1$ and $\n_p=\det$ on
$A_p=\M_2(\QQ_p)$,
\begin{align*}
[\OOO_p:\Lambda_p]&=
\big|\ZZ_p/(p^{2\max\{a,0\}}\,\ZZ_p)\big|\;
\big|\ZZ_p/(p^{\max\{a,0\}+\max\{b,0\}}\,\ZZ_p)\big|^2\;
\big|\ZZ_p/(p^{2\max\{b,0\}}\,\ZZ_p)\big|
\\ & =p^{4(\max\{a,0\}+\max\{b,0\})}=p^{-4(\min\{-a,0\}+\min\{-b,0\})}
=p^{-4\,\nu_p(\n_p(\OOO_p z^{-1}+\OOO_p))}\;,
\end{align*}
as wanted.

Now assume that $p$ divides $D_A$, so that $A_p$ is a division
algebra. Let $\nu=\nu_p\circ\n_p$, which is a discrete valuation on
$A_p$, whose valuation ring is $\OOO_p$ (see for instance \cite[page
34]{Vigneras80}). The left ideals of $\OOO_p$ are two-sided ideals.
Let $\pi$ be a uniformizer of $\OOO_p$. Note that the residual field
$\OOO_p/\pi\OOO_p$ has order $p^2$, and that $\n_p(\OOO_p)=1$ and
$\n_p(\pi)=p$. We have
$$
\Lambda_p=\OOO_p\;\pi^{2\max\{\nu(z),\,0\}}\;\;\;{\rm and}\;\;\;
\OOO_p\, z^{-1}+\OOO_p=
\OOO_p\;\pi^{\min\{\nu(z^{-1}),\,0\}}\;.
$$
Hence $[\OOO_p:\Lambda_p]=p^{4\max\{\nu(z),\,0\}}$ and
$\nu_p(\n_p(\OOO_p\, z^{-1}+\OOO_p))=-\max\{\nu(z),0\}$, which
is also as wanted.  \cqfd

\section{Background on Hamilton-Bianchi groups}
\label{sect:hambiagroup}

The Dieudonné determinant (see \cite{Dieudonne71,Aslaksen96}) $\Det$
is the group morphism from the group $\operatorname{GL}_2(\HH)$ of
invertible $2\times 2$ matrices with coefficients in $\HH$ to
$\RR^*_+$, defined by
\begin{equation}\label{eq:detdieud}
\Det\big(\Big(\begin{array}{cc}a& b\\c& d\end{array}\Big)\big)^2\;=
\n(a\,d)+ \n(b\,c) - \tr(a\,\ov c\,d\,\ov b) =\left\{\begin{array}{cl}
\n(ad - aca^{-1}b) &{\rm if}\; a \neq 0\\
\n(cb - cac^{-1}d) &{\rm if}\; c \neq 0\\
\n(cb - db^{-1}ab) &{\rm if}\; b \neq 0\;. \end{array}\right.
\end{equation}
It is invariant under the adjoint map $g\mapsto g^*$, by the
properties of $\n$ and $\tr$.  We will
denote by $\SLH$ the group of $2\times 2$ matrices with coefficients
in $\HH$ with Dieudonn\'e determinant $1$, which equals the group of
elements of (reduced) norm $1$ in the central simple algebra
$\M_2(\HH)$ over $\RR$, see \cite[Sect.~9a]{Reiner75}). We refer for
instance to \cite{Kellerhals03} for more information on $\SLH$.

The group $\SLH$ acts linearly on the left on the right $\HH$-module
$\HH\times\HH$. Let $\PP^1_r(\HH)= (\HH\times\HH-\{0\})/
\HH^\times$ be the right projective line of $\HH$, identified as usual
with the Alexandrov compactification $\HH\cup\{\infty\}$ where $[1:0]
=\infty$ and $[x:y]=xy^{-1}$ if $y\neq 0$. The projective action of
$\SLH$ on $\PP^1_r(\HH)$, induced by its linear action on $\HH\times
\HH$, is then the action by homographies on $\HH\cup\{\infty\}$
defined by
$$
\Big(\begin{array}{cc} a & b \\ c & d\end{array}\Big)\cdot z =
\left\{\begin{array}{ll} (az+b)(cz+d)^{-1} &
{\rm if}\; z\neq \infty,-c^{-1}d \\
ac^{-1} & {\rm if}\; z=\infty, c\neq 0\\
\infty & {\rm otherwise~.}\end{array}\right.
$$
This action by homographies induces a faithful left action of
$\PSLH=\SLH/\{\pm \id\}$ on $\HH\cup\{\infty\}$.

\medskip The group $\PSLH$ is very useful to study $5$-dimensional
real hyperbolic geometry, for the following reason. Let us endow $\HH$
with its usual Euclidean metric $ds^2_\HH$ (invariant under
translations, with $(1,i,j,k)$ orthonormal). We will denote by
$x=(z,r)$ a generic point in $\HH\times\,]0,+\infty[$, and by
$r:x\mapsto r(x)$ the second projection in this product.  For the real
hyperbolic space $\HH^5_\RR$ of dimension $5$, we will use the upper
halfspace model $\HH\times\,]0,+\infty[$ with Riemannian metric
$ds^2(x)= \frac{ds^2_\HH(z)+dr^2}{r^2}$ at the point $x=(z,r)$, whose
volume form is
\begin{equation}\label{eq:formevolume}
\dvol_{\HH^5_\RR}(x)=\frac{\dvol_\HH(z) \;dr}{r^5}\;.
\end{equation}
The space at infinity $\partial_\infty\hcr$ is hence
$\HH\cup\{\infty\}$.

By the Poincaré extension procedure (see for instance \cite[Lemma
6.6]{ParPau10GT}), the action of $\SLH$ by homographies on
$\partial_\infty\hcr$ extends to a left action on $\hcr$ by
\begin{equation}\label{eq:Poincareextension}
 \Big(\begin{array}{cc} a & b \\ c
    & d\end{array}\Big)\cdot (z,r)= 
\Big(\;\frac{(az+b)\,\overline{(cz+d)}+a\,\overline{c}\,r^2}
{\n(cz+d)+r^2\n(c)}, \, \frac r{\n(cz+d)+r^2\n(c)}\,\Big)\;.
\end{equation}
In this way, the group $\PSLH$ is identified with the group of
orientation preserving isometries of $\hcr$. Note that the isomorphism
$\PSLH\simeq \operatorname{SO}_0(1,5)$ is one of the isomorphisms
between connected simple real Lie groups of small dimensions in
E.~Cartan's classification.

\medskip Given an order $\OOO$ in a definite quaternion algebra $A$
over $\QQ$, define the {\it Hamilton-Bianchi group} as $\Ga_\OOO=\SLO=
\SLH \cap\M_2(\OOO)$. Note that since the norm $\n$ takes integral
values on $\OOO$, and since the Dieudonn\'e determinant is a group
morphism, we have $\GLO=\SLO$. The Hamilton-Bianchi group $\Ga_\OOO$
is a (nonuniform) arithmetic lattice in the connected real Lie group
$\SLH$ (see for instance \cite[page 1104]{ParPau10GT} for details). In
particular, the quotient real hyperbolic orbifold $\Ga_\OOO\bs\hcr$
has finite volume. The action by homographies of $\Ga_\OOO$ preserves
the right projective space $\PP^1_r(\OOO) = A\cup\{\infty\}$, which is
the set of fixed points of the parabolic elements of $\Ga_\OOO$ acting
on $\hcr\cup\partial_\infty \hcr$.

\medskip
\noindent\rem 
For every $(u,v)$ in $\OOO\times\OOO-\{(0,0)\}$, consider the two
left ideals of $\OOO$
$$
I_{u,v}=\OOO u+\OOO v\;,\;\;K_{u,v}=\Big\{\begin{array}{cl}
\OOO u\cap\OOO v &{\rm if~} uv\neq 0\;,\\ 
\OOO&{\rm otherwise.}\end{array}
$$
The map
$$
\Ga_\OOO\bs\PP^1_r(\OOO)\ra (_\OOO\!\I\times \;_\OOO\!\I)\;,
$$
which associates, to the orbit of $[u:v]$ in $\PP^1_r(\OOO)$ under
$\Ga_\OOO$, the couple of ideal classes $([I_{u,v}],[K_{u,v}])$ is a
bijection.  To see this, let $\ell_{u,v}:\OOO\times\OOO\ra \OOO$ be
the morphism of left $\OOO$-modules defined by $(o_1,o_2)\mapsto
o_1u+o_2v$.  The map $w\mapsto (wu^{-1},-wv^{-1})$ is an isomorphism
of left $\OOO$-modules from $\OOO u\cap\OOO v$ to the kernel of
$\ell_{u,v}$ if $uv\neq 0$. The result then follows for instance from
\cite[Satz 2.1, 2.2]{KraOse90}), which says that the map $[u:v]
\mapsto ([\operatorname{im}\ell_{u,v}],[\ker \ell_{u,v}])$ induces a
bijection from $\Ga_\OOO\bs\PP^1_r(\OOO)$ into $_\OOO\!\I\times
\;_\OOO\!\I$.

In particular, the number of cusps of $\Ga_\OOO$ (or the number of
ends of $\Ga_\OOO\bs\hcr$) is the square of the class number $h_A$
of $A$.

\section{Background on  binary Hamiltonian 
forms}
\label{sect:indefbinhamform}

With $V$ the right $\HH$-module $\HH\times\HH$, a {\it binary
  Hamiltonian form} $f: V \ra \RR$ is a map $X\mapsto \phi(X,X)$ where
$\phi:V\times V\ra \HH$ is a Hermitian form on $V$ with the
conjugation as the anti-involution of the ring $\HH$.  That is,
$\phi(X\lambda,Y)= \overline{\lambda}\;\phi(X,Y)$, $\phi(X+X',Y)=
\phi(X,Y)+\phi(X',Y)$, $\phi(Y,X)= \overline{\phi(X,Y)}$ for
$X,X',Y\in V$ and $\lambda\in\HH$. Our convention of sesquilinearity
on the left is the opposite of Bourbaki's unfortunate one in
\cite{Bourbaki59}. Equivalently, a binary Hamiltonian form $f$ is a
map $\HH \times \HH \ra \RR$ with
$$
f(u, v) = a\n(u) + \tr(\ov u\, b\, v) + c\n(v)\,
$$
whose {\it coefficients} $a=a(f)$ and $c=c(f)$ are real, and $b=b(f)$
lies in $\HH$. Note that $f((u,v)\lambda)=\n(\lambda)f(u,v)$. The {\it
  matrix} $M(f)$ of $f$ is the Hermitian matrix
$\Big(\begin{array}{cc}a& b\\ \ov b& c\end{array}\Big)$, so that
$f(u,v)=\Big(\begin{array}{c}\!u\!\\ \!v\! \end{array}\Big)^*\;
\Big(\begin{array}{cc}a& b\\\ov b& c\end{array}\Big)\;
\Big(\begin{array}{c} \!u\!\\ \!v\!\end{array}\Big)$.  The {\it
  discriminant} of $f$ is 
$$\Delta=\Delta (f) = \n(b)- ac.$$ Note that
the sign convention of the discriminant varies in  the
references. An easy computation shows that the Dieudonné determinant
of $M(f)$ is equal to $|\Delta |$. If $a\neq 0$, then
\begin{equation}\label{eq:factorformbinham}
f(u, v) = a\Big(\n\big(u+\frac{bv}{a}\big) 
-\frac{\Delta}{a^2}\n(v)\Big)\;.
\end{equation}
Hence the form $f$ is {\it indefinite} (that is, $f$ takes both
positive and negative values) if and only if $\Delta$ is positive, and
$\Delta$ is then equal to the Dieudonné determinant of $M(f)$.  By
Equation \eqref{eq:factorformbinham}, the form $f$ is {\it positive
  definite} (that is, $f(x)\geq 0$ with equality if and only if $x=0$)
if and only if $a>0$ and $\Delta<0$.

The linear action on the left on $\HH\times\HH$ of the group $\SLH$
induces an action on the right on the set of binary Hermitian forms
$f$ by precomposition, that is, by $f\mapsto f\circ g$ for every $g\in
\SLH$. The matrix of $f\circ g$ is $M(f\circ g)= g^*\,M(f)\,g$.  Since
the Dieudonné determinant is a group morphism, invariant under the
adjoint map (and since $f\circ g$ is indefinite if and only if $f$
is), we have, for every $g\in \SLH$, 
\begin{equation}\label{eq:invardiscrim}
\Delta(f\circ g)=\Delta(f)\;.
\end{equation}

Given an order $\OOO$ in a definite quaternion algebra over $\QQ$, a
binary Hamiltonian form $f$ is {\it integral} over $\OOO$ if its
coefficients belong to $\OOO$.  Note that such a form $f$ takes
integral values on $\OOO\times \OOO$. The lattice $\Ga_\OOO=\SLO$ of
$\SLH$ preserves the set of indefinite binary Hamiltonian forms $f$
that are integral over $\OOO$. The stabilizer in $\Ga_\OOO$ of such a
form $f$ is its {\it group of automorphs}
$$
\autom=\{g\in\Ga_\OOO\;:\;f\circ g=f\}\;.
$$

For every indefinite binary Hamiltonian form $f$, with $a=a(f)$,
$b=b(f)$ and $\Delta=\Delta(f)$, let
$$
\C_\infty(f)=\{[u:v]\in\PP^1_r(\HH)\;:\;f(u,v)=0\}\;\;\;{\rm and}
$$
$$
\C(f)=\{(z,r)\in\HH\times\,]0,+\infty[\;:\;f(z,1)+a\,r^2=0\}\;.
$$
In $\PP^1_r(\HH)=\HH\cup\{\infty\}$, the set $\C_\infty(f)$ is the
$3$-sphere of center $-\frac{b}{a}$ and radius
$\frac{\sqrt{\Delta}}{|a|}$ if $a\neq 0$, and it is the union of
$\{\infty\}$ with the real hyperplane $\{z\in\HH\;:\;
\tr(\overline{z}b)+c=0\}$ of $\HH$ otherwise. The map $f\mapsto
\C_\infty(f)$ induces a bijection between the set of indefinite binary
Hamiltonian forms up to multiplication by a nonzero real factor and
the set of $3$-spheres and real hyperplanes in $\HH\cup\{\infty\}$.
The action of $\SLH$ by homographies on $\HH\cup\{\infty\}$ preserves
this set of $3$-spheres and real hyperplanes, and the map $f\mapsto
\C_\infty(f)$ is (anti-)equivariant for the two actions of $\SLH$, in
the sense that, for every $g\in\SLH$,
\begin{equation}\label{eq:antiequiv}
\C_\infty(f\circ g)= g^{-1}\,\C_\infty(f)\;.
\end{equation}

\medskip 
Given a finite index subgroup $G$ of $\SLO$, an integral binary
Hamiltonian form $f$ is called {\it $G$-reciprocal} if there exists an
element $g$ in $G$ such that $f\circ g=-f$. We define $R_G(f)=2$ if
$f$ is $G$-reciprocal, and $R_G(f)=1$ otherwise.  The values of $f$
are positive on one of the two components of $\PP^1_r(\HH)-
\C_\infty(f)$ and negative on the other. As the signs are switched by
precomposition by an element $g$ as above, the $G$-reciprocity of the
form $f$ is equivalent to saying that there exists an element of $G$
preserving $\C_\infty(f)$ and exchanging the two complementary
components of $\C_\infty(f)$.

\section{Using Eisenstein series to compute hyperbolic 
volumes}
\label{sect:Eisenstein}

Let $\OOO$ be a maximal order in a definite quaternion algebra $A$
over $\QQ$.

\medskip In this section, we compute $\Vol(\PSLO\bs \hcr)$ using a
method which goes back in dimension $2$ to Rankin-Selberg's method
\cite{Rankin39,Selberg40} of integrating Eisenstein series on
fundamental domains and ``unfolding'', generalized by
\cite{Langlands66b} to the lattice of $\ZZ$-points of any connected
split semi-simple algebraic group over $\QQ$.  We follow the approach
of \cite[page 261-262]{Sarnak83} in dimension 3. We refer to the
Appendix for a completely different proof by V.~Emery of the same
result.

\btheo\label{theo:maincomputvol}
Let $\OOO$ be a maximal order in a definite quaternion
algebra $A$ over $\QQ$ with discriminant $D_A$. Then
$$
\Vol(\PSLO\bs \hcr)=
\frac{\zeta(3)\prod_{p|D_A}(p^3-1)(p-1)\,}{11520}\;.
$$
\etheo

\dem It is well known (see for instance \cite[Sect.~6.3,
Ex. (3)]{ParPau10GT}) that there exists $\underline{G}$, a connected
semi-simple linear algebraic group over $\QQ$, such that
$\underline{G}(\RR)=\SLH$, $\underline{G}(\QQ) =\SL_2(A)$ and
$\underline{G}(\ZZ)=\SLO$. Let $\underline{P}$ be the parabolic
subgroup of $\underline{G}$, defined over $\QQ$, such that
$\underline{P}(\RR)$ is the upper triangular subgroup of $\SLH$.  By
Borel's finiteness theorem (see for instance \cite{Borel66}), the set
$\SLO\bs\SL_2(A)/\underline{P}(\QQ)$ is finite, and we will fix a
subset $\Rep$ in $\SL_2(A)$ which is a system of representatives of
this set of double cosets.

Let $\Ga=\SLO$. For every $\alpha\in \Rep$, let
$\Ga_\alpha=\underline{P}(\RR)\cap (\alpha^{-1} \Ga\alpha)$ and let
$\Ga'_\alpha$ be its subgroup of unipotent elements. The group
$\alpha\Ga_\alpha\alpha^{-1}$ is the stabilizer of the parabolic fixed
point $\alpha\infty$ in $\Ga$. The action of $\Ga_\alpha$ on
$\HH\cup\{\infty\}$ by homographies preserves $\infty$ and is
cocompact on $\HH$. If $\alpha^{-1}=\Big(\!\begin{array}{cc}a& b\\c&
  d\end{array}\!\Big)$ and $\alpha=\Big(\!\begin{array}{cc}\wt a& \wt
  b\\\wt c& \wt d\end{array}\!\Big)$, let $\uuu_\alpha=c\OOO+d\OOO$,
which is a right fractional ideal of $\OOO$, and $\vvv_\alpha=\OOO\wt
a+\OOO\wt c$, which is a left fractional ideal of $\OOO$.

\medskip For every $\alpha\in \Rep$, the {\it Eisenstein series} of
the arithmetic group $\Ga$ for the cusp at infinity $\alpha\infty$ is
the map $E_\alpha:\hcr\times \;]4,+\infty[ \; \ra \RR$ defined
by
$$
E_\alpha(x,s)=
\sum_{\ga\in (\alpha\Ga_\alpha\alpha^{-1})\bs\Ga} r(\alpha^{-1}\ga x)^s\;.
$$
The summation does not depend on the choice of representatives of the
left cosets in $(\alpha\Ga_\alpha\alpha^{-1})\bs\Ga$ since
$\Ga_\alpha$ preserves $\infty$ and the Euclidean height $r$.  The
{\em Eisenstein series} of $\OOO$ is (for $x=(z,r)\in\hcr$ and
$s\in\CC$ with $\Re \;s>4$)
$$
\wh E(x,s)=\sum_{(c,\,d)\in\OOO\times\OOO-\{0\}} 
\Big(\frac{r}{\n(cz+d)+r^2\n(c)}\Big)^s\;.
$$

The following result is a concatenation of results proven in
\cite{KraOse90}.

\btheo [Krafft-Osenberg] \label{theo:rappelsKO} (i) The Eisenstein
series $E_\alpha(x,s)$, $\alpha\in \Rep$, and $\wh E(x,s)$ converge
absolutely and uniformly on compact subsets of $\{s\in\CC\;:\;\Re\;
s>4\}$, uniformly on compact subsets of $x\in \hcr$. They are
invariant by the action of $\Ga$ on the first variable.

(ii) The map $s\mapsto \wh E(x,s)$ admits a meromorphic extension
to $\CC$, having only one pole, which is at $s=4$ and is simple with
residue
\begin{equation}\label{eq:residuEhat}
\Res_{s=4} \wh E(x,s)=\frac{8\,\pi^4}{3\,{D_A}^2}\;. 
\end{equation}  
Furthermore, if $c(\alpha,s)= \n(\uuu_\alpha)^s\;
\zeta(\uuu_\alpha^{-1},\frac{s}{2})$ for every $\alpha\in \Rep$, then
\begin{equation}\label{eq:relatEhatEinftygene}
\wh E(x,s)= \sum_{\alpha\in \Rep} \;c(\alpha,s) \;E_\alpha(x,s)
\;. 
\end{equation}

(iii) For every $\alpha,\beta\in \Rep$, there exist a map $s\mapsto
\varphi_{\alpha,\beta}(s)$ with $(s-4)\varphi_{\alpha,\beta} (s)$
bounded for $s>4$ near $s=4$, and a measurable map $(x,s)\mapsto
\Phi_{\alpha,\beta} (x,s)$ such that $(s-4)\Phi_{\alpha,\beta}(x,s)$
is bounded by an integrable (for the hyperbolic volume) map,
independent on $s>4$ near $s=4$, on $x\in K\times[\epsilon,+\infty[$
where $K$ is a compact subset of $\HH$ and $\epsilon>0$, such that
$$
E_\alpha(\beta x,s)=\delta_{\alpha,\beta}\,r^s+
\varphi_{\alpha,\beta}(s)r^{4-s}+\Phi_{\alpha,\beta}(x,s)\;,
$$
with $\delta_{\alpha,\beta}=1$ if $\alpha=\beta$ and
$\delta_{\alpha,\beta} = 0$ otherwise.  
\etheo

\dem We are using Langlands' convention for the Eisenstein series,
hence with $\Ga'_{\alpha}$ the subgroup of unipotent elements of
$\Ga_{\alpha}$, our Eisenstein series $E_\alpha$ is obtained from the
one used in \cite{KraOse90} by replacing $\alpha$ by $\alpha^{-1}$ and
by multiplying by $\frac{1} {[\Ga_{\alpha}: \Ga'_{\alpha}]}$.

The part of claim (i) concerning the series  $E_\alpha(x,s)$,
$\alpha\in \Rep$, is \cite[Satz 3.2]{KraOse90}. The rest follows
from \cite[Satz~4.2]{KraOse90} with $M=\OOO$.  The claim (ii) follows
from \cite[Koro.~5.6 a)]{KraOse90} with $M=\OOO$, recalling that the
reduced discriminant of any maximal order of $A$ is equal to the
reduced discriminant of $A$. Equation \eqref{eq:relatEhatEinftygene}
follows from \cite[Satz~4.3]{KraOse90}, recalling the above changes
between our $E_\alpha$ and the one in \cite{KraOse90}.  The claim (iii)
follows from \cite[Satz~3.3]{KraOse90}, again replacing $\beta$ by
$\beta^{-1}$, and using the second equation in \cite[page
85]{MagObeSon66} to control the modified Bessel function.
\cqfd

\bigskip By a {\it fundamental domain} for a smooth action of a
countable group $G$ on a smooth manifold $N$, we mean a subset $F$ of
$N$ such that $F$ has negligible boundary, the subsets
$g\stackrel{\circ}{F}$ for $g\in G$ are pairwise disjoint, and
$N=\bigcup_{g\in G}\;gF$.

Here is a construction of a fundamental domain $\F$ for $\Ga$ acting
on $\hcr$ that will be useful in this section (and is valid for any
discrete subgroup of isometries of $\HH^n_\RR$ with finite covolume
which is not cocompact).  Let $\P$ be the set of parabolic fixed
points of $\Ga$. By the structure of the cusp neighbourhoods, there
exists a family $(\H_p)_{p\in \P}$ of pairwise disjoint closed
horoballs in $\hcr$, equivariant under $\Ga$ (that is $\ga
\H_p=\H_{\ga p}$ for every $\ga\in\Ga$), with $\H_p$ centered at
$p$. The {\it cut locus of the cusps} $\Sigma$ is the piecewise
hyperbolic polyhedral complex in $\hcr$ consisting of the set of
points outside the union of these horoballs which are equidistant to
at least two of these horoballs (it is independent of the choice of
this family when there is only one orbit of parabolic fixed points).
Each connected component of the complement of $\Sigma$ contains one
and only one of these horoballs, is at bounded Hausdorff distance of
it, is invariant under the stabilizer in $\Ga$ of its point at
infinity, and is precisely invariant under the action of $\Ga$.
Recall that a subset $A$ of a set endowed with an action of a group
$G$ is said to be {\it precisely invariant} under this group if for
every $g\in G$, if $gA\cap A$ is nonempty, then $gA=A$.

For every $\beta\in \Rep$, let $\D_\beta$ be a compact fundamental domain
for the action of $\Ga_{\beta}$ on $\HH$, let $\wt \F_\beta$ be the
closure of the component of the complement of $\Sigma$ containing
$\H_{\beta\infty}$, and define
$$
\F_\beta=\wt \F_\beta\cap \beta( \D_\beta\times\,]0,+\infty[\,)\;.
$$
Then $\F_\beta$ is a closed fundamental domain for the action of $\beta
\Ga_\beta \beta^{-1}$ on $\wt \F_\beta$, and there exists a continuous
map $\sigma'_\beta:\D_\beta\ra\;]0,+\infty[$ (which hence has a positive
lower bound) such that
\begin{equation}\label{eq:formevoiscusp}
\beta^{-1}\F_\beta=\{(z,r)\in \hcr\;:\;
z\in\D_{\beta},\;r\geq \sigma'_\beta(z)\}\;.
\end{equation}
Now define
\begin{equation}\label{eq:decompdomfond}
\F=\bigcup_{\beta\in \Rep}\F_\beta\;.
\end{equation}
Since $\Rep$ is a system of representatives of the cusps, $\F$ is a 
fundamental domain of $\Ga$ acting on $\hcr$.

Note that, for every $\alpha\in \Rep$, there exists a continuous map
$\sigma_\alpha:\D_\alpha\ra\;[0,+\infty[$ (hence with a finite upper
bound), with only finitely many zeros, such that, since $\alpha^{-1}
\F$ is a fundamental domain for the action of $\alpha^{-1} \Ga\alpha$
on $\hcr$,
\begin{equation}\label{eq:regroupdomfond}
\bigcup_{\ga\in(\alpha^{-1} \Ga\alpha-\Ga_\alpha)} \ga\alpha^{-1}\F 
=\Ga_\alpha\;\{(z,r)\in \hcr\;:\; 
z\in\D_{\alpha},\;r< \sigma_\alpha(z)\}\;.
\end{equation}

\medskip For every $\alpha\in \Rep$, let 
$$
b_\alpha(s)=
\int_{\F}\big(E_\alpha(x,s)- r(\alpha^{-1}x)^s \big)\;
\dvol_{\hcr}(x)\;.
$$
When $s>4$, we have 
\begin{align*}
b_\alpha(s)&=\int_{\F}\Big(
\sum_{\ga\in(\alpha\Ga_{\alpha}\alpha^{-1})\bs\Ga}\;r(\alpha^{-1}\ga x)^s -
r(\alpha^{-1}x)^s\Big)\; \dvol_{\HH^5_\RR}(x)\\ &=
\int_{\F}\sum_{\ga\in\Ga_{\alpha}\bs(\alpha^{-1} \Ga\alpha-\Ga_\alpha)}\;
r(\ga\alpha^{-1} x)^s\; \dvol_{\HH^5_\RR}(x)\\ &=
\sum_{\ga\in\Ga_{\alpha}\bs(\alpha^{-1} \Ga\alpha-\Ga_\alpha)}
\int_{\F}\;r(\ga\alpha^{-1} x)^s\; \dvol_{\HH^5_\RR}(x)\\ &=
\sum_{\ga\in\Ga_{\alpha}\bs(\alpha^{-1} \Ga\alpha-\Ga_\alpha)}
\int_{\ga\alpha^{-1} \F}\;r(x)^s\; \dvol_{\HH^5_\RR}(x)\\ &= 
\int_{\bigcup_{\ga\in\Ga_{\alpha}\bs(\alpha^{-1} \Ga\alpha-\Ga_\alpha)}
\ga\alpha^{-1}\F} r(x)^s \; \dvol_{\HH^5_\RR}(x)\\ &=
\int_{z\in \D_{\alpha}}\int_0^{\sigma_\alpha(z)}
r^{s-5}\;drdz=
\int_{z\in \D_{\alpha}}\frac{\sigma_\alpha(z)^{s-4}}{s-4}
\;dz\;,
\end{align*}
using for the succession of equations, respectively,
the definition of $E_\alpha$, the change of variables
$\alpha^{-1}\ga\alpha\ra\ga$, Fubini's theorem for positive
functions, the invariance of the volume under the isometric change
of variables $\ga\alpha^{-1}x\ra x$, $\sigma$-additivity, and 
the equations \eqref{eq:regroupdomfond} and \eqref{eq:formevolume} and
the invariance of the Euclidean height function $r$ under
$\Ga_\alpha$. 

For any $\alpha\in\Rep$, the map ${\sigma_\alpha}^{s-4}$ converges
pointwise, as $s\ra 4^+$, to the map on $\D_\alpha$ with value $0$ at
the finitely many points where $\sigma_\alpha$ vanishes, and with
value $1$ otherwise.  Since $\D_\alpha$ is compact and
${\sigma_\alpha}^{s-4}$ is uniformly bounded from above, Lebesgue's
dominated convergence theorem gives
$$
\lim_{s\ra 4^+}\;(s-4)\,b_\alpha(s)= 
\Vol(\D_\alpha)=\Vol(\Ga_{\alpha}\bs\HH)\;.
$$
Therefore by using Equation \eqref{eq:relatEhatEinftygene}, the map
$$
s\mapsto b(s)= \int_{\F}\big(\wh E(x,s)-\sum_{\alpha\in \Rep} 
\;c(\alpha,s)\,r(\alpha^{-1}x)^s
\big)\; \dvol_{\hcr}(x)=
\sum_{\alpha\in \Rep} 
\;c(\alpha,s)\;b_\alpha(s)\
$$
satisfies
\begin{equation}\label{eq:residue3}
\lim_{s\ra 4^+}\;(s-4)\,b(s)= \sum_{\alpha\in \Rep} 
\;c(\alpha,4)\;\Vol(\Ga_{\alpha}\bs\HH)\;,
\end{equation}
since $s\mapsto c(\alpha,s)$ is holomorphic for $\Re\;s>2$.

On the other hand, let us prove that we may permute the limit as $s\ra
4^+$ and the integral defining $(s-4)b(s)$. Using the equations
\eqref{eq:relatEhatEinftygene} and \eqref{eq:decompdomfond}, and an
isometric, hence volume preserving, change of variable, we have
\begin{align*}
  b(s)&=\sum_{\alpha,\beta\in \Rep}\;c(\alpha,s)\;\int_{\F_\beta}
  \big(E_\alpha(x,s)-r(\alpha^{-1}x)^s\big)\;\dvol_{\HH^5_\RR}(x)\\
  & =\sum_{\alpha,\beta\in \Rep}\;c(\alpha,s)\;\int_{\beta^{-1}\F_\beta}
  \big(E_\alpha(\beta x,s)-r(\alpha^{-1}\beta x)^s\big)\;
  \dvol_{\HH^5_\RR}(x)\;.
\end{align*}
If $x\in\beta^{-1}\F_\beta$, then $r(x)$ is bounded from below by a
positive constant by the construction of $\F_\beta$, hence
$r(x)^{4-s}$ is bounded from above for every $s\geq 4$.  If
$\alpha\neq \beta$ and $x\in\beta^{-1}\F_\beta$, then
$r(\alpha^{-1}\beta x)^s$ is bounded from above for every $s\geq 0$,
since $\alpha^{-1}\F_\beta$ is bounded in $\HH\times\RR$ by
construction. Hence since $\beta^{-1}\F_\beta$ has finite hyperbolic
volume, by Theorem \ref{theo:rappelsKO} (iii) separating the case
$\alpha=\beta$ and the case $\alpha\neq \beta$, by Lebesgue's
dominated convergence theorem, we may permute the limit as $s\ra 4^+$
and the integral on $\beta^{-1}\F_\beta$ for the hyperbolic volume
applied to $(s-4) \big(E_\alpha(\beta x,s)-r(\alpha^{-1}\beta
x)^s\big)$. By a finite summation, we may indeed permute the limit as
$s\ra 4^+$ and the integral defining $(s-4)b(s)$.

Therefore, by Equation \eqref{eq:residuEhat},
\begin{equation}\label{eq:residue1}
\lim_{s\ra 4^+}\;(s-4)\,b(s)=\frac{8\,\pi^4}{3\,{D_A}^2}\;
\Vol(\PSLO\backslash\HH^5_\RR)\;.
\end{equation}

Finally, since for every $\rho\in A-\{0\}$, the element
$\ga_\rho=\Big(\begin{array}{cc}\rho& -1\\1& 0\end{array}\Big)$ of
$\SL_2(A)$ maps $\infty$ to $\rho$, the element $\alpha\in\Rep$ may be
choosen to be either $\id$ or $\ga_\rho$ for some $\rho\in A$.
In the first case, $\uuu_\alpha=\OOO$ and $\Ga'_\alpha$ acts on $\HH$ as the
$\ZZ$-lattice $\OOO$, so that, by Equation \eqref{eq:volumeHHparOOO},
since the subgroup $\{\pm \id\}$ of $\Ga_\alpha$ is the kernel
of its action on $\HH$,
$$
\n(\uuu_\alpha)^4\;\Vol(\Ga_{\alpha}\bs\HH)=
\frac{2\,\Vol(\OOO\bs\HH)}{[\Ga_\alpha:\Ga'_\alpha]}=
\frac{D_A}{2\,[\Ga_\alpha:\Ga'_\alpha]}\;.
$$
In the second case, $\alpha^{-1}= \Big(\begin{array}{cc}0& 1\\-1&
  \rho\end{array}\Big)$, so that $\uuu_\alpha=\OOO\rho+\OOO$ and
$\Ga'_\alpha$ acts on $\HH$ as the $\ZZ$-lattice $\Lambda=\OOO\cap
\rho^{-1}\OOO\cap \OOO\rho^{-1}\cap \rho^{-1}\OOO\rho^{-1}$ as proved
in Lemma \ref{lem:computvolcusp}. By Lemma \ref{lem:chenevier} applied
with $z=\rho^{-1}$ and by Equation \eqref{eq:volumeHHparOOO}, we hence
have, 
$$
\n(\uuu_\alpha)^4\,\Vol(\Ga_{\alpha}\bs\HH)=
\n(\uuu_\alpha)^4\,[\OOO:\Lambda]\,
\frac{2\Vol(\OOO\bs\HH)}{[\Ga_\alpha:\Ga'_\alpha]}=
\frac{D_A}{2\,[\Ga_\alpha:\Ga'_\alpha]}\;.
$$
Therefore, by the definition of $c(\alpha,s)$,
\begin{equation}\label{eq:residue2}
\sum_{\alpha\in \Rep} 
\;c(\alpha,4)\;\Vol(\Ga_{\alpha}\bs\HH)=
\frac{D_A}{2}\;\sum_{\alpha\in \Rep}\;
\frac{\zeta(\uuu_\alpha^{-1},2)}{[\Ga_\alpha:\Ga'_\alpha]} \;.
\end{equation}
Combining the equations \eqref{eq:residue3}, \eqref{eq:residue1} and
\eqref{eq:residue2}, we have
\begin{equation}\label{eq:oldcomputvolEis}
\Vol(\PSLO\bs \hcr)=
\frac{3\,{D_A}^3}{16\,\pi^4}\sum_{\alpha\in \Rep} 
\;\frac{\zeta(\uuu_\alpha^{-1},2)}{[\Ga_\alpha:\Ga'_\alpha]}\;.
\end{equation}

\blemm\label{lem:computindex}
(1) For every $\alpha\in \Rep$, we have
$$
[\Ga_\alpha:\Ga'_\alpha]=\;\big|\OOO_r(\uuu_\alpha^{-1})^\times\big|
\;\big|\OOO_r(\vvv_\alpha)^\times\big|\;.
$$

(2) The map from $\Rep$ to $_\OOO\!\I\times
\;_\OOO\!\I$ defined by
$$
\alpha\mapsto ([\vvv_\alpha],[\uuu_\alpha^{-1}])
$$
is a bijection.
\elemm

\dem (1) Let $\Ga_\alpha^+=\big\{\ga\in\Ga_\alpha\;:\; (0 \; 1)\ga=(0 \;
1)\big\}$, and $\Ga_\alpha^-=\big\{\ga\in\Ga_\alpha\;:\;
\ga\Big(\begin{array}{c}\!\! 1\!\! \\ \!\!0\!\!\end{array}\Big)
=\Big(\begin{array}{c}\!\! 1\!\! \\ \!\!0\!\!\end{array}\Big)\big\}$,
which are normal subgroups of $\Ga_\alpha$, whose union generates
$\Ga_\alpha$.  By the top of page 434 in \cite{KraOse90} (keeping in
mind that our $\alpha$ is the inverse of the $\alpha$ in
\cite{KraOse90}), we have $[\Ga_\alpha^+: \Ga'_\alpha]
=\big|\OOO_r(\vvv_\alpha)^\times\big|$.  Similarly, $[\Ga_\alpha^-:
\Ga'_\alpha]= \big| \OOO_\ell(\uuu_\alpha)^\times \big|$.  Note that
$\Ga'_\alpha$ is a normal subgroup of $\Ga_\alpha^-,\Ga_\alpha^+$ and
$\Ga_\alpha$, such that $\Ga_\alpha^-\cap\Ga_\alpha^+ = \Ga_\alpha'$.
Hence the product map from $\Ga_\alpha^-\times\Ga_\alpha^+$ to
$\Ga_\alpha'$ induces a bijection from $(\Ga_\alpha^-/\Ga_\alpha')
\times(\Ga_\alpha'\bs\Ga_\alpha^+)$ to $\Ga_\alpha/\Ga_\alpha'$, since
$\Ga_\alpha/\Ga_\alpha'$ is abelian. In particular,
$[\Ga_\alpha:\Ga'_\alpha]=\;\big|\OOO_\ell(\uuu_\alpha)^\times\big|
\;\big|\OOO_r(\vvv_\alpha)^\times\big|$. Using Equation
\eqref{eq:gauchedroite}, the result follows.  

\medskip (2) Since these matrices act transitively on $A$ by
homographies, we may assume that every $\alpha\in\Rep$ either is the
identity element $\id$, or has the form $\Big(\begin{array}{cc}\!\!
  \rho_\alpha & -1\!\! \\ \!\!1 & 0\!\!\end{array}\Big)$ for some
$\rho_\alpha\in A^\times$.  Then $\alpha^{-1}$ is either $\id$ or
$\Big(\begin{array}{cc}\!\!  0 & 1\!\! \\ \!\!-1 &
  \rho_\alpha\!\!\end{array}\Big)$. So that
$\uuu_\alpha=\OOO+\rho_\alpha\OOO$ and
$\vvv_\alpha=\OOO\rho_\alpha+\OOO$, unless $\alpha=\id$, in which case
$\uuu_\alpha=\vvv_\alpha=\OOO$.  Since $\SL_2(A)$ acts (on the left)
transitively by homographies on $\PP^1_r(\OOO)$ with stabilizer of
$[1:0]$ equal to $\underline{P}(\QQ)$, the map from $\Rep$ to
$\Ga_\OOO\bs\PP^1_r(\OOO)$ defined by $\alpha\mapsto
\Ga_\OOO\alpha[1:0]$ is a bijection.  Note that
$\alpha[1:0]=[\rho_\alpha:1]$ if $\alpha\neq \id$.  Using the notation
of the last remark of Section \ref{sect:hambiagroup}, if $\alpha\neq
\id$, we have $\vvv_\alpha =I_{\rho_\alpha,1}$ and
$[K_{\rho_\alpha,1}] = [\OOO\rho_\alpha\cap\OOO] =
[\OOO\cap\OOO\rho_\alpha^{-1}]= [\uuu_\alpha^{-1}]$ by Equation
\eqref{eq:invsomeginter}. The second assertion of this lemma then
follows from the last remark of Section \ref{sect:hambiagroup}.  \cqfd

\medskip Now, using respectively Equation \eqref{eq:KraOse436}, Lemma
\ref{lem:computindex} (1), Lemma \ref{lem:computindex} (2), the
separation of variables and Equation \eqref{eq:Deuring128}, Equation
\eqref{eq:Deuring134}, and Equation \eqref{eq:valeurzetaooo} since
$\zeta(4)=\frac{\pi^4}{90}$, we have
\begin{align}
\sum_{\alpha\in\Rep}\;
\frac{\zeta(\uuu_\alpha^{-1},2)}{[\Ga_\alpha:\Ga'_\alpha]}
&=\sum_{\alpha\in\Rep}\;
\frac{\big|\OOO_r(\uuu_\alpha^{-1})^\times\big|\;
\zeta_{[\uuu_\alpha^{-1}]}(2)}{[\Ga_\alpha:\Ga'_\alpha]} =\sum_{\alpha\in\Rep}\; 
\frac{\zeta_{[\uuu_\alpha^{-1}]}(2)}
{\big|\OOO_r(\vvv_\alpha)^\times\big|}
\nonumber
\\ & =\sum_{([I],[J])\,\in\; _\OOO\!\I\times\;_\OOO\!\I}\; 
\frac{\zeta_{[J]}(2)}
{\big|\OOO_r(I)^\times\big|}=\zeta_{A}(2)\;\sum_{[I]\,\in\; _\OOO\!\I}\; 
\frac{1}
{\big|\OOO_r(I)^\times\big|}
\nonumber\\ &=
\frac{\zeta_{A}(2)}{24}\;\prod_{p|D_A} (p-1)
=\frac{\zeta(3)\;\pi^4\;\prod_{p|D_A}(1-p^{-3})(p-1)}{2160}
\label{eq:sommepourrie}\;.
\end{align}

Theorem \ref{theo:maincomputvol} follows from the equations
\eqref{eq:oldcomputvolEis} and \eqref{eq:sommepourrie}.  
\cqfd

\bcoro \label{coro:volhycomput} Let $A$ be a definite quaternion
algebra over $\QQ$ with reduced discriminant $D_A$ and class number
$1$, and let $\OOO$ be a maximal order in $A$. Then the hyperbolic
volume of $\PSLO\bs \hcr$ is equal to
$$
\Vol(\PSLO\bs \hcr)= \frac{(D_A^3-1)(D_A-1)\,\zeta(3)}{11520}\;.
$$
\ecoro

This is an immediate consequence of Theorem
\ref{theo:maincomputvol}. But here is a proof directly from Equation
\eqref{eq:oldcomputvolEis} which avoids using the technical Lemma
\ref{lem:computindex} and the technical computation
\eqref{eq:sommepourrie}.

\medskip \dem Since the number of cusps of $\SLO$ is the square of the
class number $h_A$ of $A$ (see the remark at the end of Section
\ref{sect:hambiagroup}), the set $\Rep$ has only one element, and we
may choose $\Rep=\{\id\}$.

By definition of the Dieudonn\'e determinant and since every element of
$\OOO^\times$ has norm $1$, the stabilizer $\Ga_\infty$ of $\infty$ in
$\SLO$ is
$$
\Ga_\infty=\Big\{\Big(\begin{array}{cc}a& b\\0& d\end{array}\Big)\;:\;
a,d\in\OOO^\times,\;b\in\OOO\Big\}\;.
$$
The index in $\Ga_\infty$ of its unipotent subgroup is hence
$|\OOO^\times|^2$. By the equations \eqref{eq:valeurpartzetaideal} and
\eqref{eq:valeurzetaooo}, Corollary \ref{coro:volhycomput} follows
from Equation \eqref{eq:oldcomputvolEis}, since
$\zeta(4)=\frac{\pi^4}{90}$, since $|\OOO^\times|= \frac{24}{D_A-1}$
as seen in Equation \eqref{eq:computnombreunit}, and since $D_A$ is
prime when $h_A=1$.  \cqfd

\medskip\noindent {\bf Example. } Let $A$ be Hamilton's quaternion
algebra over $\QQ$, which satisfies $D_A=2$ and $h_A=1$. Let $\OOO$ be
Hurwitz's maximal order in $A$. Applying Corollary
\ref{coro:volhycomput}, we get
$$
\Vol(\PSLO\backslash \hcr)= \frac{7\,\zeta(3)}{11520}\;,
$$ 
exactly four times the minimal volume of a cusped hyperbolic
$5$-orbifold, as we should because the Hurwitz modular group
$\PSL_2(\OOO)$ is a subgroup of index $4$ in the group of the minimal
volume cusped hyperbolic orbifold of dimension $5$, see \cite[page
209]{Hild07} and \cite[page 186]{JohWei99}.

\section{Representing integers by binary 
Hamiltonian forms}
\label{sect:rephamil}

Let $A$ be a definite quaternion algebra over $\QQ$, and let $\OOO$ be
a maximal order in $A$.

Let us introduce the general counting function we will study.  For
every indefinite integral binary Hamiltonian form $f$ over $\OOO$, for
every finite index subgroup $G$ of $\SLO$, for every $x,y$ in $\OOO$
not both zero, and for every $s>0$, let
$$
\psi_{f,G,x,y}(s)=\card\;\;_{\mbox{$\autom\cap G$}}\bs
\big\{(u,v)\in G(x,y)\;:\;
\n(\OOO x+\OOO y)^{-1}|f(u,v)|\leq s\big\}\;.
$$
The counting function $\psi_{f,G,x,y}$ depends (besides on $f,G$) only
on the $G$-orbit of $[x:y]$ in $\PP^1_r(\OOO)$. 

Here is the notation for the statement of our main result which
follows.  Given $(x,y)\in \OOO\times \OOO$, let $\Ga_{\OOO,x,y}$ and
$G_{x,y}$ be the stabilizers of $(x,y)$ for the left linear actions of
$\Ga_\OOO= \SLO$ and $G$ respectively, and let $\uuu_{xy^{-1}}$ be the
right fractional ideal $\OOO$ if $y=0$ and $\OOO+xy^{-1}\OOO$
otherwise.  Let $\iota_G=1$ if $-\id\in G$, and $\iota_G=2$
otherwise. Note that the image of $\autom\cap G$ in $\PSLH$ is again
an arithmetic group.

\btheo\label{theo:mainversionG} Let $f$ be an integral indefinite
binary Hamiltonian form of discriminant $\Delta(f)$ over a maximal order
$\OOO$ of a definite quaternion algebra $A$ over $\QQ$. Let $x$ and
$y$ be elements in $\OOO$ not both zero, and let $G$ be a finite index
subgroup of $\Ga_\OOO=\SLO$. Then, as $s$ tends to $+\infty$, we have
the equivalence
$$
\psi_{f,G,x,y}(s)\sim \;\;
\frac{540\;\iota_G\;[\Ga_{\OOO,x,y}:G_{x,y}]\;\covol(\autom\cap G)}
{\pi^2\;\zeta(3)\;|\OOO_\ell(\uuu_{xy^{-1}})^\times|\;\Delta(f)^2
\;[\Ga_\OOO:G]\prod_{p|D_A}(p^3-1)(1-p^{-1})} \;\;\;s^4\;.
$$
\etheo

\dem Let us first recall a geometric result from \cite{ParPau} that
will be used to prove this theorem.  

Let $n\geq 2$ and let $\hnr$ be the upper halfspace model of the real
hyperbolic space of dimension $n$, with (constant) sectional curvature
$- 1$. Let $F$ be a finite covolume discrete group of isometries of
$\hnr$. Let $1\leq k\leq n-1$ and let $\C$ be a real hyperbolic
subspace of dimension $k$ of $\hnr$, whose stabilizer $F_\C$ in $F$
has finite covolume. Let $\H$ be a horoball in $\hnr$, which is
precisely invariant under $F$, with stabilizer $F_\H$.

For every $\alpha,\beta\in F$, denote by $\delta_{\alpha,\beta}$ the
common perpendicular geodesic arc between $\alpha \C$ and the
horosphere $\beta\partial\H$, and let $\ell(\delta_{\alpha,\beta})$ be
its length, counted positively if $\delta_{\alpha,\beta}$ exits
$\beta\H$ at its endpoint on $\beta\partial\H$, and negatively
otherwise. By convention, $\ell(\delta_{\alpha,\beta})=-\infty$ if the
boundary at infinity of $\alpha \C$ contains the point at infinity of
$\beta\H$.  Also define the multiplicity of $\delta_{\alpha,\beta}$ as
$m(\alpha, \beta) = 1/{\card(\alpha\,F_\C\alpha^{-1}\; \cap\; \beta\,
  F_\H\beta^{-1})}$. Its denominator is finite, if the boundary at
infinity of $\alpha \C$ does not contain the point at infinity of
$\beta\H$, since then the subgroup $\alpha\,F_\C\alpha^{-1}\; \cap\;
\beta\, F_\H\beta^{-1}$ that preserves both $\beta\H$ and $\alpha \C$,
consists of elliptic elements. In particular, there are only finitely
many elements $[g]\in F_\C\bs F/F_\H$ such that $m(g^{-1},\id)$ is
different from $1$, or equivalently such that $g^{-1}F_\C g\; \cap\;
F_\H\neq \{1\}$.  For every $t\geq 0$, define $\N(t)= \N_{F,\C,\H}(t)$
as the number, counted with multiplicity, of the orbits under $F$ in
the set of the common perpendicular arcs $\delta_{\alpha,\beta}$ for
$\alpha,\beta\in F$ with length at most $t$:
$$
\N(t)=\N_{F,\C,\H}(t)=\sum_{(\alpha,\,\beta)\,\in\, 
F\backslash ((F/F_\C)\times (F/F_\H))\;:\;
\ell(\delta_{\alpha,\beta})\leq t} m(\alpha,\beta)\;.
$$

For every $m\in\NN$, denoting by  $\SSS_m$  the unit
sphere of the Euclidean space $\RR^{m+1}$, endowed with its induced
Riemannian metric, we have the following result:

\btheo [{\cite[Coro.~4.9]{ParPau}}]\label{theo:mainBPP}
As $t\to+\infty$, we have
$$
\N(t)\sim\frac{\Vol(\SSS_{n-k-1})\Vol(F_\H\bs\H)\Vol(F_\C\bs\C)}
{\Vol(\SSS_{n-1})\Vol(F\bs\hnr)}\;e^{(n-1)t}\;.\;\;\Box
$$
\etheo

\bigskip Now, let $A,\OOO,f,G,x$ and $y$ be as in the statement of
Theorem \ref{theo:mainversionG}. We write $f$ as in Equation
\eqref{eq:form}, and denote its discriminant by $\Delta$. In order to
apply Theorem \ref{theo:mainBPP}, we first define the various objects
$n$, $k$, $F$, $\H$, and $\C$ that appear in its statement.

Let $n=5$ and $k=4$, so that $\Vol(\SSS_{n-1})=\frac{8\pi^2}{3}$ and
$\Vol(\SSS_{n-k-1}) = 2$. We use the description of $\hcr$ given in
Section \ref{sect:hambiagroup}.

For any subgroup $S$ of $\SLH$, we denote by $\ov S$ its image in
$\PSLH$, except that the image of $\autom$ is denoted by $\Pautom$. We
will apply Theorem \ref{theo:mainBPP} to $F=\ov G$.

\medskip Note that $\Vol(\,\ov G\,\bs\hcr)=[\,\overline{\Ga_\OOO}:\ov
G\,]\Vol(\,\overline{\Ga_\OOO}\,\bs\hcr)$ and
$[\,\overline{\Ga_\OOO}:\ov G\,]= \frac{1}{\iota_G}[\Ga_\OOO:G]$ by
the definition of $\iota_G$. Thus, using Theorem
\ref{theo:maincomputvol} (or Theorem \ref{theo:appendix} in the
Appendix), we have
\begin{align}
\Vol(\,\ov G\,\bs\hcr)&=
\frac{1}{\iota_G}[\Ga_\OOO:G]\;\Vol(\,\overline{\Ga_\OOO}\,\bs\hcr)=
\frac{1}{\iota_G}[\Ga_\OOO:G]\;\covol(\Ga_\OOO)\nonumber
\\ \label{eq:volhypmani}& =
\frac{\zeta(3)\;[\Ga_\OOO:G]}{11520\;\iota_G}\;\prod_{p|D_A}(p^3-1)(p-1)
\;.
\end{align}

The point $\rho=xy^{-1}\in A\cup\{\infty\}\subset\phcr$ is a parabolic
fixed point of $\overline{\Ga_\OOO}$ hence of $\ov G$.  Let
$\tau\in\;]0,1]$ and $\H$ be the horoball in $\hcr$ centered at
$\rho$, with Euclidean height $\tau$ if $y\neq 0$, and consisting of
the points of Euclidean height at least $\frac{1}{\tau}$ otherwise.
Assume that $\tau$ is small enough so that $\H$ is precisely invariant
under $\overline{\Ga_\OOO}$ hence under $\ov G$.  Such a $\tau$
exists, as seen in the construction of the fundamental domain in
Section \ref{sect:Eisenstein}. The stabilizer
$\overline{\Ga_{\OOO,\rho}}$ in $\overline{\Ga_\OOO}$ of the point at
infinity $\rho$ is equal to the stabilizer
$(\,\overline{\Ga_\OOO}\,)_\H$ of the horoball $\H$.

\rem If $\rho=\infty$ and $G=\Ga_\OOO$, we may take $\tau=1$ by
\cite[Prop.~5]{Kellerhals03}. Then by an easy hyperbolic geometry
computation, since the index in $(\,\overline{\Ga_\OOO}\,)_\H$ of the
subgroup of translations by elements of $\OOO$ is
$\frac{|\OOO^\times|^2}{2}$, and by using Equation
\eqref{eq:volumeHHparOOO}, we have
$$
\Vol((\,\overline{\Ga_\OOO}\,)_\H\bs\H)= 
\frac{1}{4}\Vol((\,\overline{\Ga_\OOO}\,)_\H\bs\partial\H)=
\frac{1}{2|\OOO^\times|^2}\Vol(\OOO\bs\HH)=\frac{D_A}{8|\OOO^\times|^2}\;.
$$
The following lemma will allow us to generalize this formula.

\blemm \label{lem:computvolcusp} Let $\Lambda'_{\OOO,\rho}=\OOO\cap
\rho^{-1}\OOO\cap \OOO\rho^{-1}\cap \rho^{-1}\OOO\rho^{-1}$ if
$x,y\neq 0$, and $\Lambda'_{\OOO,\rho}=\OOO$ otherwise. Then
$\Lambda'_{\OOO,\rho}$ is a $\ZZ$-lattice in $\HH$ and we have
\begin{equation}\label{eq:volcusp}
\Vol(\,\ov G_\H\bs\H)= 
\frac{\tau^4\,[(\,\overline{\Ga_\OOO}\,)_\H:\ov G_\H]}
{4\,|\OOO_\ell(\uuu_\rho)^\times|\,
[(\,\overline{\Ga_\OOO}\,)_\H:\overline{\Ga_{\OOO,x,y}}]}
\;\Vol(\Lambda'_{\OOO,\rho}\bs\HH)\;.
\end{equation}
\elemm

\dem If $y= 0$, let $\ga_\rho=\id$, otherwise let 
$$
\ga_\rho=\Big(\!\begin{array}{cc}\rho& -1\\1&
  \;\,0\end{array}\!\Big)\in\SLH\;.
$$
Note that $\ga_\rho^{-1}$ maps $\rho$ to $\infty$ and $\H$ to the
horoball $\H_\infty$ consisting of the points in $\hcr$ with
Euclidean height at least $\frac{1}{\tau}$.

Let $\ga=\Big(\!\begin{array}{cc}a& b\\c& d\end{array}\!\Big)$ and
$\ga'= \Big(\!\begin{array}{cc}1& b'\\0& 1\end{array}\!\Big)$ be in
$\SLH$.  If $y=0$, we have $\ga_\rho^{-1}\ga\ga_\rho=\ga'$ if and only
if $a=1$, $b=b'$, $c=0$, $d=1$. If $y\neq 0$, by an easy
computation, we have $\ga_\rho^{-1}\ga\ga_\rho=\ga'$ (that is
$\ga\ga_\rho= \ga_\rho \ga'$) if and only if
\begin{equation}\label{eq:tasdegalite} 
c=-b'\;,\;\;a=1-\rho b'\;,\;\;d=1+b'\rho\;,\;\;
b=\rho b'\rho\;.
\end{equation}
In particular, if $x,y\neq 0$, if $\ga\in \SLO$ and $\ga'=
\ga_\rho^{-1} \ga \ga_\rho\in \SL_2(A)$ is unipotent upper triangular,
then these equations imply respectively that $b'$ belongs to $\OOO$,
$\rho^{-1}\OOO$, $\OOO\rho^{-1}$ and $ \rho^{-1}\OOO\rho^{-1}$,
therefore $b'\in\Lambda'_{\OOO,\rho}$. If $x=0$ or $y=0$, we also have
$b'\in\OOO=\Lambda'_{\OOO,\rho}$

Conversely, if $b'\in \Lambda'_{\OOO,\rho}$, then define $a,b,c,d$ by
the equations \eqref{eq:tasdegalite} if $y\neq 0$, and
by $a=1$, $b=b'$, $c=0$, $d=1$ otherwise, so that $a,b,c,d\in\OOO$.
Let $\ga=\Big(\!\begin{array}{cc}a& b\\c& d\end{array}\!\Big)$. If
$y\neq 0$, note that if $c=0$, then $\ga=\id$ and otherwise $cb
-cac^{-1}d=-1$, so that $\ga\in\SLO$ by Equation \eqref{eq:detdieud}.
If $y\neq 0$, the equations \eqref{eq:tasdegalite} imply that
$\ga_\rho^{-1}\ga \ga_\rho$ is a unipotent upper triangular element of
$\SLO$, and this is also the case if $y=0$.

The abelian group $\Lambda'_{\OOO,\rho}$ is a $\ZZ$-lattice in $\HH$,
as an intersection of at most four $\ZZ$-lattices in $A$.  Since an
isometry preserves the volume for the first equality, by an easy
hyperbolic volume computation for the second one, and by the previous
computation of the unipotent upper triangular subgroup
$\Ga'_{\ga_\rho}$ of $\ga_\rho^{-1}\Ga_{\OOO,x,y} \ga_\rho$ for the
last one, we have
\begin{align*}
  \Vol(\,\overline{\Ga_{\OOO,x,y}}\,\bs\H)&=
  \Vol\big((\ga_\rho^{-1}\,\overline{\Ga_{\OOO,x,y}}\,\ga_\rho\,)
  \bs\H_\infty\big)=\frac{1}{4}\;\Vol\big(
  (\ga_\rho^{-1}\,\overline{\Ga_{\OOO,x,y}}\,\ga_\rho\,)
  \bs\partial\H_\infty\big)\\ & = \frac{\tau^4}{4}\;\Vol\big(
  (\ga_\rho^{-1}\,\overline{\Ga_{\OOO,x,y}}\,\ga_\rho\,)
  \bs\HH\big)=\frac{\tau^4}{4\,[\ga_\rho^{-1}\Ga_{\OOO,x,y}
    \ga_\rho:\Ga'_{\ga_\rho}]}\;\Vol( \Lambda'_{\OOO,\rho}\bs\HH)\;.
\end{align*}
With the notation of the proof of Lemma 10(1), we have
$\ga_\rho^{-1}\Ga_{\OOO,x,y} \ga_\rho=\Ga^-_{\ga_\rho}$, hence
$[\ga_\rho^{-1}\Ga_{\OOO,x,y} \ga_\rho:\Ga'_{\ga_\rho}]=
|\OOO_\ell(\uuu_\rho)^\times|$. Since covering arguments
yield
\begin{align*}
\Vol(\,\ov G_\H\bs\H)&=[(\,\overline{\Ga_\OOO}\,)_\H:\ov G_\H]\;
\Vol\big((\,\overline{\Ga_\OOO}\,)_\H\bs\H\big)=
\frac{[(\,\overline{\Ga_\OOO}\,)_\H:\ov G_\H]}
{[(\,\overline{\Ga_\OOO}\,)_\H:\overline{\Ga_{\OOO,x,y}}\,]}\;
\Vol(\,\overline{\Ga_{\OOO,x,y}}\,\bs\H)\;,
\end{align*}
the result follows. \cqfd

\bigskip Let us resume the proof of Theorem \ref{theo:mainversionG}.
Let $\C=\C(f)$, which is indeed a real hyperbolic hyperplane in
$\hcr$, whose set of points at infinity is $\C_\infty(f)$ (hence
$\infty$ is a point at infinity of $\C(f)$ if and only if the first
coefficient $a=a(f)$ of $f$ is $0$). Note that $\C$ is invariant under
the group $\autom$ by Equation \eqref{eq:antiequiv} (which implies
that $\C(f\circ g)= g^{-1}\C(f)$ for every $g\in\SLO$). The arithmetic
group $\autom$ acts with finite covolume on $\C(f)$, its finite
subgroup $\{\pm\id\}$ acting trivially. By definition,
$$
\covol(\autom\cap G)= \Vol\big(\;_{\mbox{$\Pautom\cap \ov G$}}
\bs\C(f)\big)\;.
$$  
Note that $\covol(\autom\cap G)$ depends only on the $G$-orbit of $f$,
by Equation \eqref{eq:antiequiv} and since $\operatorname{SU}_{f \circ
  g} (\OOO)=g^{-1}\autom g$ for every $g\in \SLO$.  By its definition,
$R_G(f)$ is the index of the subgroup $\Pautom\cap \ov G$ in $\ov
G_\C$, hence
\begin{equation}\label{eq:voltotgeod}
\Vol(\,\ov G_\C\bs\C)=\frac{1}{R_G(f)}\covol(\autom\cap G)\;.
\end{equation}

\medskip The last step of the proof of Theorem \ref{theo:mainversionG}
consists in relating the two counting functions $\psi_{f,G,x,y}$ and
$\N_{\ov G,\C,\H}$, in order to apply Theorem \ref{theo:mainBPP}.  

For every $g\in\SLH$, let us compute the hyperbolic length of the
common perpendicular geodesic arc $\delta_{g^{-1},\id}$ between the
real hyperbolic hyperplane $g^{-1}\C$ and the horoball $\H$, assuming
that they do not meet.  We use the notation $\ga_\rho,\H_\infty$
introduced in the proof of Lemma \ref{lem:computvolcusp}. Since
$\ga_\rho^{-1}$ sends the horoball $\H$ to the horoball $\H_\infty$,
it sends the common perpendicular geodesic arc between $g^{-1}\C$ and
$\H$ to the (vertical) common perpendicular geodesic arc between
$\ga_\rho^{-1}g^{-1}\C$ and $\H_\infty$. Let $r$ be the Euclidean
radius of the $3$-sphere $\C_\infty(f\circ g\circ \ga_\rho)$, which is
the image by $\ga_\rho^{-1}$ of the boundary at infinity of $g^{-1}\C$
by Equation \eqref{eq:antiequiv}. Denoting by $a(f\circ g\circ
\ga_\rho)$ the coefficient of $\n(u)$ in $f\circ g\circ
\ga_\rho(u,v)$, we have, by Equation \eqref{eq:invardiscrim},
$$
r=\frac{\sqrt{\Delta(f\circ g\circ
    \ga_\rho)}} {|a(f\circ g\circ \ga_\rho)|}=
\frac {\sqrt{\Delta}}{|f\circ
g\circ \ga_\rho(1,0)|}=\frac {\sqrt{\Delta}}{|f\circ
g(\rho,1)|}
=\frac{\n(y)\;\sqrt{\Delta}}{|f\circ g(x,y)|}\;,
$$
if $y\neq 0$ and $r=\frac{\n(x)\;\sqrt{\Delta}}{|f\circ g(x,y)|}$
otherwise. An immediate computation gives
\begin{equation}\label{eq:calcullongarcperpcomm}
\ell(\delta_{g^{-1},\id})= \ell(\ga_\rho^{-1}\delta_{g^{-1},\id})= 
\ln\frac{1}{\tau}-\ln r
=\ln\frac{|f\circ g(x,y)|}{\tau\,\n(y)\;\sqrt{\Delta}}\;,
\end{equation}
if $y\neq 0$ and $\ell(\delta_{g^{-1},\id})=\ln\frac{|f\circ
  g(x,y)|}{\tau\,\n(x)\;\sqrt{\Delta}}$ otherwise. With the
conventions that we have taken, these formulas are also valid if
$g^{-1}\C$ and $\H$ meet.

\medskip Recall that there are only finitely many elements $[g]\in
\overline{G}_\C \bs \overline{G}/\overline{G}_\H$ such that $g^{-1}
\overline{G}_\C \,g\; \cap\; \overline{G}_\H$ is different from
$\{1\}$ or such that the multiplicity $m(g^{-1},\id)$ is different
from $1$. If $y\neq 0$, using Equation
\eqref{eq:calcullongarcperpcomm} for the second line below, Lemma 7 in
\cite{ParPau11CJM} for the third one, and Theorem \ref{theo:mainBPP}
applied to $F=\ov G$ for the fifth one, we hence have, as $s$ tends to
$+\infty$,
$$
\begin{aligned}
  \psi_{f,G,x,y}(s) &=\card\;\big\{[g]\in(\autom\cap G)\bs G/
  G_{x,y}\;:\;
\n(\OOO x+\OOO y)^{-1}  |f\circ g(x,y)|\leq s\big\}\\
  &= \card\;\big\{[g]\in(\Pautom\cap\ov G)\bs \ov G/\,\overline
  {G_{x,y}}\;:\;\ell(\delta_{g^{-1},\id})\leq 
\ln{\textstyle
\frac{s\;\n(\OOO x+\OOO y)}{\tau\,\n(y)\;\sqrt{\Delta}}\big\}}
\\ &\sim R_G(f)\;[\,\ov G_\H:\overline{G_{x,y}}\,]\;
  \card\;\big\{[g]\in \ov G_\C\bs \ov G/\,\ov G_\H\;:\;
  \ell(\delta_{g^{-1},\id})\leq 
\ln{\textstyle\frac{s\;\n(\OOO \rho+\OOO )}{\tau\,\sqrt{\Delta}}}
\big\}\\ &
\sim R_G(f)\;[\,\ov G_\H:\overline{G_{x,y}}\,]\;\N_{\ov G,\C,\H}\big(
\ln{\textstyle\frac{s\;\n(\OOO \rho+\OOO)}{\tau\,\sqrt{\Delta}}}
\big)\\ &
\sim R_G(f)\;[\,\ov G_\H:\overline{G_{x,y}}\,]\;
\frac{6\Vol(\,\ov G_\H\bs\H)\Vol(\,\ov G_\C\bs\C)}
{8\pi^2\Vol(\,\ov G\bs\hnr)}\;
\Big(\frac{s\;\n(\OOO \rho+\OOO)}{\tau\,\sqrt{\Delta}}\Big)^4\;.
\end{aligned}
$$
We replace the three volumes in the above computation by their
expressions given in the equations \eqref{eq:volhypmani},
\eqref{eq:volcusp}, \eqref{eq:voltotgeod}. We simplify the obtained
expression using the following two remarks. Firstly,
$$
[\,\ov G_\H:\overline{G_{x,y}}\,]\;\frac{
[(\,\overline{\Ga_\OOO}\,)_\H:\,\ov G_\H]}
{[(\,\overline{\Ga_\OOO}\,)_\H:\,\overline{\Ga_{\OOO,x,y}}]}=
\frac{[(\,\overline{\Ga_\OOO}\,)_\H:\overline{G_{x,y}}\,]}
{[(\,\overline{\Ga_\OOO}\,)_\H:\,\overline{\Ga_{\OOO,x,y}}]}
=[\,\overline{\Ga_{\OOO,x,y}}:\overline{G_{x,y}}\,]=
[\Ga_{\OOO,x,y}:G_{x,y}\,]\;.
$$
Secondly, we claim that 
\begin{equation}\label{eq:claimfincount}
\Vol(\Lambda'_{\OOO,\rho}\bs\HH)\n(\OOO\rho+\OOO)^4=
\frac{D_A}{4}\;.
\end{equation}
If $x=0$, then $\Lambda'_{\OOO,\rho}=\OOO$, hence this claim is true,
by Equation \eqref{eq:volumeHHparOOO} and since $\n(\OOO)=1$.  Otherwise,
Claim \eqref{eq:claimfincount} follows from Lemma \ref{lem:chenevier}
with $z=\rho^{-1}$, since, by the definition of
$\Lambda'_{\OOO,\rho}$,
\begin{align*}
\Vol(\Lambda'_{\OOO,\rho}\bs\HH)
\n(\OOO\rho+\OOO)^4&= 
\Vol(\Lambda\bs\HH)\n(\OOO z^{-1}+\OOO)^4= 
\Vol(\OOO\bs\HH)[\OOO:\Lambda]\n(\OOO z^{-1}+\OOO)^4\;,
\end{align*}
and by Equation \eqref{eq:volumeHHparOOO}.

\medskip This concludes the proof of Theorem \ref{theo:mainversionG}
if $y\neq 0$. The case $y=0$ is similar to the case $x=0$. \cqfd

\medskip Let us give a few corollaries of Theorem
\ref{theo:mainversionG}. The first one below follows by taking
$G=\SLO$ in Theorem \ref{theo:mainversionG}.

\bcoro\label{coro:GegalSLO} Let $f$ be an integral indefinite binary
Hamiltonian form of discriminant $\Delta(f)$ over a maximal order
$\OOO$ of a definite quaternion algebra $A$ over $\QQ$. Let $x$ and
$y$ be elements in $\OOO$ not both zero. Then, as $s$ tends to
$+\infty$, we have the equivalence
$$
\psi_{f,\SLO,x,y}(s)\sim \;\; \frac{540\;\covol(\autom)}
{\pi^2\;\zeta(3)\;|\OOO_\ell(\uuu_{xy^{-1}})^\times|\;
  \Delta(f)^2\;\prod_{p|D_A}(p^3-1)(1-p^{-1})} \;\;\;s^4\;. \;\;\;\Box
$$
\ecoro

\noindent{\bf Remarks.}
(1) Recall that by the last remark of Section \ref{sect:hambiagroup},
the map from $\SLO\bs\PP^1_r(\OOO)$ to $_\OOO\!\I\times \;_\OOO\!\I$
which associates, to the orbit of $[u:v]$ in $\PP^1_r(\OOO)$ under
$\SLO$, the couple of ideal classes $([I_{u,v}],[K_{u,v}])$ is a
bijection. The counting function $\psi_{f,\SLO,x,y}$ hence depends
only on $([I_{x,y}],[K_{x,y}])$.

\medskip (2) Given two left fractional ideals $\mmm$ and $\mmm'$ of
$\OOO$, let
 $$
\psi_{f,\mmm,\mmm'}(s)=\card\;\;_{\mbox{$\autom$}}\bs
\big\{(u,v)\in \mmm\times\mmm\;:\;
\frac{|f(u,v)|}{\n(\mmm)}\leq s\;,\;\;I_{u,v}=\mmm\;,\;\;
[K_{u,v}]=[\mmm'] \ \big\}\;.
$$
Note that this counting function depends only on the ideal classes of
$\mmm$ and $\mmm'$.

\bcoro\label{coro:doublefracideal} Let $f$ be an integral indefinite
binary Hamiltonian form of discriminant $\Delta(f)$ over a maximal
order $\OOO$ of a definite quaternion algebra $A$ over $\QQ$. Let
$\mmm$ and $\mmm'$ be two left fractional ideals in $\OOO$. Then as
$s$ tends to $+\infty$, we have the equivalence
$$
\psi_{f,\mmm,\mmm'}(s)\sim \;\;
\frac{540\;\covol(\autom)}
{\pi^2\;\zeta(3)\;|\OOO_r(\mmm')^\times|\;\Delta(f)^2\;
\prod_{p|D_A}(p^3-1)(1-p^{-1})}\;\;\;s^4\;. \;\;\;\Box
$$
\ecoro

\dem By Remark (1), we have
$$
\psi_{f,\mmm,\mmm'}=\psi_{f,\SLO,x,y},
$$
where $(x,y)$ is any nonzero element of $\OOO\times\OOO$ such that
$[I_{x,y}]=[\mmm]$ and $[K_{x,y}]=[\mmm']$. By Equation
\eqref{eq:gauchedroite} and \eqref{eq:invsomeginter}, if $xy\neq 0$, we
have
$$
|\OOO_\ell(\uuu_{xy^{-1}})^\times|=|\OOO_r({\uuu_{xy^{-1}}}^{-1})^\times|=
|\OOO_r(\OOO \cap\OOO yx^{-1})^\times|=|\OOO_r(K_{x,y})^\times|\,.
$$
The first and last terms are also equal if $xy= 0$. Hence the result
follows from Corollary \ref{coro:GegalSLO}.  \cqfd

\medskip (3) With $\psi_{f,\mmm}$ the counting function
defined in the Introduction, we have
\begin{equation}\label{eq:somdeunadeuidealfrac}
\psi_{f,\mmm}=\sum_{[\mmm']\;\in\;_\OOO\!\I}\;\psi_{f,\mmm,\mmm'}\;.
\end{equation} 
Therefore, since
$\sum_{[\mmm']\;\in\;_\OOO\!\I}
\frac{1}{|\OOO_r(\mmm')^\times|}=
  \frac{1}{24}\;\prod_{p|D_A} (p-1)$ by Equation
\eqref{eq:Deuring134}, Theorem \ref{theo:mainintro} in the
Introduction follows from Corollary \ref{coro:doublefracideal}.

\medskip
We say $u,v\in\OOO\times\OOO$ are {\it relatively prime} if one of the
following equivalent (by Remark (1) above) conditions is satisfied~:

(i) there exists $g\in \SLO$ such that $g(1,0)=(u,v)$;

(ii) there exists $u',v'\in\OOO$ such that 
$\n(uv')+\n(u'v)-\tr(u\,\ov v \,v'\,\overline{u'})=1$;

(iii) the $\OOO$-modules $I_{u,v}$ and $K_{u,v}$ are isomorphic (as
$\OOO$-modules) to $\OOO$.

\noindent 
We denote by $\P_{\OOO}$ the set of couples of relatively prime
elements of $\OOO$. 

\bcoro\label{coro:caspartrelprime} Let $f$ be an integral indefinite
binary Hamiltonian form over a maximal order $\OOO$ in a definite
quaternion algebra $A$ over $\QQ$, and let $G$ be a finite index
subgroup of $\Ga_\OOO=\SLO$.  Then, as $s$ tends to $+\infty$, we have
the equivalence
$$
\card\;\;_{\mbox{$\autom\cap G$}}\bs
\big\{(u,v)\in \P_\OOO\;:\;
|f(u,v)|\leq s\big\}
$$
$$
\sim \;\; \frac{540\;\iota_G\;
[\Ga_{\OOO,1,0}:G_{1,0}]\,\covol(\autom\cap G)}
{\pi^2\;\zeta(3)\;|\OOO^\times|\;\Delta(f)^2\;[\Ga_\OOO:G]\;
\prod_{d|D_A}(p^3-1)(1-p^{-1})}\;\;\;s^4\;.
$$
\ecoro

\dem This follows from Theorem \ref{theo:mainversionG} by taking $x=1$
and $y=0$.  
\cqfd

\medskip Let us now give a proof of Corollary \ref{coro:introspunooo}
in the Introduction.

\medskip
\noindent {\bf Proof of Corollary \ref{coro:introspunooo}.} Consider
the integral indefinite binary Hamiltonian form $f$ over $\OOO$
defined by $f(u,v)=\tr(\ov u\,v)$, with matrix $M(f)=
\Big(\begin{array}{cc}0& 1\\1& 0\end{array}\Big)$ and discriminant
$\Delta(f)=1$.  Its group of automorphs is
$$
\SpO=\Big\{g\in\SLO\;:\;
g^*\,\Big(\begin{array}{cc}0& 1\\1&
  0\end{array}\Big)\,g=\Big(\begin{array}{cc}0& 1\\1&
  0\end{array}\Big)\Big\}\;,
$$
which is an arithmetic lattice in the symplectic group over the
quaternions $\operatorname{Sp}_1(\HH)$. We have
$$
\C(f)=\{(z,r)\in\HH\times\,]0,+\infty[\;:\;\tr (z)= 0\}\;.
$$

The hyperbolic volume of the quotient of $\{(z,r)\in\HH\times\,
]0,+\infty[ \;:\;\tr (z)= 0\}$ by $\SpO$ has been computed as the main
result of \cite{BreHelm96}, yielding
$$
\covol(\SpO)=\frac{\pi^2}{1080}\prod_{p|D_A}(p^2+1)(p-1)\;,
$$
where $p$ ranges over the primes dividing $D_A$. 

Corollary \ref{coro:introspunooo} in the Introduction then follows from
Theorem \ref{theo:mainintro}  with $\mmm=\OOO$. \cqfd

\medskip \rem Theorem \ref{theo:mainversionG} and its Corollary
\ref{coro:caspartrelprime} allow the asymptotic study of the counting
of representations satisfying congruence properties. For instance, let
$\I$ be a (nonzero) twosided ideal in an order $\OOO$ in a definite
quaternion algebra $A$ over $\QQ$. Let $\Ga_\I$ be the kernel of the
map $\SLO\ra\operatorname{GL}_2(\OOO/\I)$ of reduction modulo $\I$ of
the coefficients, and $\Ga_{\I,0}$ the preimage of the upper
triangular subgroup by this map. Then applying Corollary
\ref{coro:caspartrelprime} with respectively $G=\Ga_\I$ or
$G=\Ga_{\I,0}$, we get an asymptotic equivalence as $s\ra+\infty$ of
the number of relatively prime representations $(u,v)$ of integers
with absolute value at most $s$ by a given integral binary Hamiltonian
form, satisfying the additional congruence properties $\{u\equiv 1
\mod \I$, $v\equiv 0\mod \I\}$ or $\{v\equiv 0\mod \I\}$. To give an
even more precise result, the computation of the indices of $\Ga_\I$
and $\Ga_{\I,0}$ in $\SLO$ would be needed.

\section{Geometric reduction theory of binary 
Hamiltonian forms}
\label{sect:geomredtheo}

Let $\OOO$ be a (not necessarily maximal) order in a definite
quaternion algebra $A$ over $\QQ$.

Let $\Q$ be the $6$-dimensional real vector space of binary
Hamiltonian forms, $\Q^+$ the open cone of positive definite ones,
$\Q^{\pm}$ the open cone of indefinite ones, $\Q(\OOO)$ the discrete
subset of the ones that are integral over $\OOO$, and
$$
\Q^+(\OOO)=\Q^+\cap \Q(\OOO), \;\;\;\;\Q^\pm(\OOO)=\Q^\pm\cap
\Q(\OOO)\;.
$$ 
For every $\Delta\in\ZZ-\{0\}$, let $\Q(\Delta)=\{f\in \Q\;:\;
\Delta(f)=\Delta\}$, $\Q(\OOO,\Delta)=\Q(\Delta)\cap \Q(\OOO)$ and
$$
\Q^+(\OOO,\Delta)=\Q(\Delta)\cap \Q^+(\OOO), \;\;
\Q^\pm(\OOO,\Delta)=\Q(\Delta)\cap\Q^\pm(\OOO)\;.
$$

The group $\RR^*_+$ acts on $\Q^+$ by multiplication; we will denote
by $[f]$ the orbit of $f$ and by $\overline{\Q}^+$ the quotient space
$\Q^+/\RR^*_+$. Similarly, the group $\RR^*$ acts on $\Q^\pm$ by
multiplication; we will denote by $[f]$ the orbit of $f$ and by
$\overline{\Q}^\pm$ the quotient space $\Q^\pm/\RR^*$.  The right
action of $\SLH$ on $\Q$ preserves $\Q(\Delta)$, $\Q^+$ and
$\Q^{\pm}$, commuting with the actions of $\RR^*_+$ and $\RR^*$ on
these last two spaces. The subgroup $\SLO$ preserves $\Q(\OOO)$,
$\Q^+(\OOO)$, $\Q^\pm(\OOO)$, $\Q^+(\OOO,\Delta)$,
$\Q^\pm(\OOO,\Delta)$.

Let $\C(\hcr)$ be the space of totally geodesic hyperplanes of $\hcr$,
with the Hausdorff distance on compact subsets.

\bprop\label{prop:equivar}
(1) The map $\Phi:\overline{\Q}^+\ra \hcr$ defined by
$$
[f]\mapsto  \Big(\,-\frac{b(f)}{a(f)},
\frac{\sqrt{-\Delta(f)}}{a(f)}\;\Big)
$$
is a homeomorphism, which is (anti-)equivariant for the actions of
$\SLH$~: For every $g\in\SLH$, we have $\Phi([f\circ g])=
g^{-1}\Phi([f])$.

(2) The map $\Psi:\overline{\Q}^\pm\ra \C(\hcr)$ defined by
$$
[f]\mapsto  \C(f)
$$
is a homeomorphism, which is (anti-)equivariant for the actions of
$\SLH$~: For every $g\in\SLH$, we have $\Psi([f\circ g])
=g^{-1}\Psi([f])$.  
\eprop

Note that $\Phi([f])$ may be geometrically understood as the pair of
the center and the imaginary radius of the imaginary sphere with
equation $f(z,1)=0$, that is $\n(z+\frac{b(f)}{a(f)})=
-\Big(\frac{\sqrt{-\Delta(f)}}{a(f)}\Big)^2$.

\medskip
\dem (1) Since $a=a(f)>0$ and $\Delta=\Delta(f)<0$ when $f$ is a
positive definite binary Hamiltonian form, the map $\Phi$ is well
defined and continuous.  Since the orbit by $\RR^*_+$ of a positive
definite binary Hamiltonian form has a unique element $f$ such that
$a(f)=1$, and since $c(f)$ then is equal to $n(b(f))-\Delta$, the map
$\Phi$ is a bijection with continuous inverse $(z,r)\mapsto [f_{z,r}]$
where $f_{z,r}:(u,v)\mapsto \n(u)-\tr(\overline{u}zv) +(\n(z)+r^2)\n(v)$.

To prove the equivariance property of $\Phi$, we could use Equation
\eqref{eq:Poincareextension} and the formula for the inverse of an
element of $\SLO$ given for instance in \cite{Kellerhals03}, but the
computations are quite technical and even longer than below. Hence we
prefer to use the following lemma to decompose the computations.
%
%
%
%

\blemm \label{lem:rootgroup}
The group (even the monoid) $\SLH$ is generated by the elements
$\Big(\!\begin{array}{cc}0& -1\\1& \;\,0\end{array}\!\Big)$ and
$\Big(\!\begin{array}{cc}1& \beta\\0& 1\end{array}\!\Big)$ with
$\beta\in\HH$.  
\elemm

This is a consequence of a general fact about connected semisimple
real Lie groups and their root groups, but the proof is short
(and is one way to prove that the Dieudonné determinant of
$\Big(\!\begin{array}{cc}\alpha& \beta\\\ga& \delta\end{array}\!\Big)$
is $\n(\ga \beta-\ga \alpha \ga^{-1}\delta)$ if $\ga\neq 0$).

\medskip \dem This follows from the following facts, where
$\alpha,\beta,\ga,\delta\in\HH$. If $\alpha\neq 0$, then
$$
\Big(\!\begin{array}{cc}\alpha& \beta\\0& \delta\end{array}\!\Big)=
\Big(\!\begin{array}{cc}\alpha& 0\\0& \delta\end{array}\!\Big)
\Big(\!\begin{array}{cc}1& \alpha^{-1}\beta\\0& 1\end{array}\!\Big)\;,
\;\;\;{\rm and}
$$
$$
\Big(\!\begin{array}{cc}\alpha& 0\\0& \alpha^{-1}\end{array}\!\Big)=
\Big(\!\begin{array}{cc}1 & -\alpha\\ 0& 1\end{array}\!\Big)
\Big(\!\begin{array}{cc}0& -1\\1& \;\,0\end{array}\!\Big)
\Big(\!\begin{array}{cc}1 & -\alpha^{-1}\\ 0& 1\end{array}\!\Big)
\Big(\!\begin{array}{cc}0& -1\\1& \;\,0\end{array}\!\Big)
\Big(\!\begin{array}{cc}1 & -\alpha\\ 0& 1\end{array}\!\Big)
\Big(\!\begin{array}{cc}0& -1\\1& \;\,0\end{array}\!\Big).
$$
If $\n(\alpha\delta)=1$, there there exists $u,v\in\HH^\times$ such that
$\alpha\delta=uvu^{-1}v^{-1}$, and
$$
\Big(\!\begin{array}{cc}\alpha& 0\\0& \delta\end{array}\!\Big)=
\Big(\!\begin{array}{cc}u& 0\\0& u^{-1}\end{array}\!\Big)
\Big(\!\begin{array}{cc}v& 0\\0& v^{-1}\end{array}\!\Big)
\Big(\!\begin{array}{cc}(vu)^{-1}& 0\\0& vu\end{array}\!\Big)
\Big(\!\begin{array}{cc}\delta^{-1}& 0\\0& \delta\end{array}\!\Big) \;.
$$
If $\ga \neq 0$, then 
$$
\Big(\!\begin{array}{cc}\alpha& \beta\\\ga & \delta\end{array}\!\Big)=
\Big(\!\begin{array}{cc}1& \alpha \ga^{-1}\\0& 1\end{array}\!\Big)
\Big(\!\begin{array}{cc}0& -1\\1& 0\end{array}\!\Big)
\Big(\!\begin{array}{cc}\ga & 0\\0& 
-\beta+\alpha \ga^{-1}\delta\end{array}\!\Big)
\Big(\!\begin{array}{cc}1& \ga^{-1}\delta\\0& 1\end{array}\!\Big)
\;.\;\;\;\Box
$$

\medskip Now, to prove the equivariance property, one only has to
prove it for the elements of the generating set of $\SLH$ given in the
above lemma. Given $f\in\Q^+$, let $M=\Big(\!\begin{array}{cc}a&
  b\\\overline{b} & c\end{array}\!\Big)$ be the matrix of $f$ and
$\Delta=\Delta(f)$. Note that the matrix of $f\circ g$ is $g^*Mg$.

If $g=\Big(\!\begin{array}{cc}1& \beta\\0& 1 \end{array} \!\Big)$, we
have $a(f\circ g) =a$ and $b(f\circ g)=a\beta+b$.  Since
$$
g^{-1}\cdot \Big(\!-\frac{b}{a}, \frac{\sqrt{-\Delta}}{a}\Big)= 
\Big (\!-\frac{b}{a} -\beta, \frac{\sqrt{-\Delta}}{a}\Big)=
\Big(\!-\frac{b(f\circ g)}{a(f\circ g)}, \frac{\sqrt{-\Delta(f\circ
    g)}}{a(f\circ g)}\Big)
$$ 
by Equation \eqref{eq:invardiscrim}, the result follows in this case.

%
If $g=\Big(\!\begin{array}{cc}0& -1\\1& 0\end{array} \!\Big)$,
then $a(f\circ g)=c$ and $b(f\circ g)= -\,\overline{ b}$.  By Equation
\eqref{eq:Poincareextension}, for every $(z,h)\in\hcr$, we have
$g^{-1}\cdot (z,r)= (\frac{-\overline{z}}{\n(z)+r^2},
\frac{r}{\n(z)+r^2})$. Therefore, since $\Delta=\n(b)-ac$,
$$
g^{-1}\cdot \Big(\!-\frac{b}{a}, \frac{\sqrt{-\Delta}}{a}\Big)
=\Big(\frac{-(-\frac{\overline{b}}{a})}{\frac{\n(b)}{a^2}+
\frac{-\Delta}{a^2}},
\frac{\frac{\sqrt{-\Delta}}{a}}{\frac{c}{a}}\Big)=
\Big(\!-\frac{b(f\circ g)}
{a(f\circ g)}, \frac{\sqrt{-\Delta(f\circ g)}}{a(f\circ g)}\Big)\;.
$$
The equivariance property of $\Phi$ follows.

\medskip (2) We have already seen that $\Psi$ is a bijection. Its
equivariance property follows from Equation \eqref{eq:antiequiv}. Let
$a=a(f),b=b(f),c=c(f)$ and $\Delta=\Delta(f)$. Since $\C(f)=\{(z,r)\in
\hcr\;:\;\n(az+b)+a^2r^2=\Delta\}$ if $a\neq 0$ and $\C(f)=\{(z,r)\in
\hcr\;:\;\tr(\overline{z}b)+c=0\}$ otherwise, the map $\Psi$ is
clearly a homeomorphism. \cqfd

\bigskip In order to define a geometric notion of reduced binary
Hamiltonian form, much less is needed than an actual fundamental
domain for the group $\SLO$ acting on $\hcr$. Though it might increase
the number of reduced elements, this will make the verification that a
given binary form is reduced much easier (see the end of this
Section). Indeed, due to the higher dimension, the number of
inequalities is much larger than the one for $\SL_2(\ZZ)$ or for
$\SL_2(\OOO_K)$ where $\OOO_K$ is the ring of integers of an imaginary
quadratic number field $K$ (see for instance
\cite{Zagier81,BucVol07,ElsGruMen98}).

\medskip 
For $n\geq 2$, let us denote by $\|z\|$ the usual Euclidean norm on
$\RR^{n-1}$. Consider the upper halfspace model of the real hyperbolic
$n$-space $\HH^n_\RR$, whose underlying manifold is $\RR^{n-1}\times
\;]0,+\infty[$, so that $\partial_\infty\HH^n_\RR=\RR^{n-1}\cup
\{\infty\}$.  A {\it weak fundamental domain} for the action of a
finite covolume discrete subgroup $\Ga$ of isometries of $\HH^n_\RR$
is a subset $\F$ of $\HH^n_\RR$ such that
\begin{enumerate}
\item[(i)] $\bigcup_{g\in\Ga}\;g\F=\HH^n_\RR$ ;
\item[(ii)] there exists a compact subset $K$ in $\RR^{n-1}$ such that $\F$
  is contained in $K\times\;]0,+\infty[$ ;
\item[(iii)] there exists $\kappa,\epsilon>0$ and a finite set $Z$ of
  parabolic fixed points of $\Ga$ such that $\F=\{(z,r)\in\F\;:\;r\geq
  \epsilon\}\cup\bigcup_{s\in Z}\E_s$ where
  $\E_s\subset\{(z,r)\in\F\;:\;\|z-s\|\leq \kappa \,r^2\}$.
\end{enumerate}
Note that a weak fundamental domain for a finite index subgroup of
$\Ga$ is a weak fundamental domain for $\Ga$.

When $\infty$ is a parabolic fixed point of $\Ga$, an example of a weak
fundamental domain is any Ford fundamental domain of $\Ga$, whose
definition we now recall.

Given any isometry $g$ of $\HH^n_\RR$ such that $g\infty\neq \infty$,
the {\it isometric sphere} of $g$ is the $(n-2)$-sphere $S_g$ of
$\RR^{n-1}$ which consists of the points at which the tangent map of
$g$ is a Euclidean isometry.  We then define $S_g^+$ as the set of
points in $\HH^n_\RR$ that are in the closure of the unbounded
component of the complement of the hyperbolic hyperplane whose
boundary is $S_g$. For instance, if $g=\Big(\!\begin{array}{cc}\alpha&
  \beta\\\ga& \delta\end{array} \!\Big)\in\SLH$, then $g\infty\neq
\infty$ if and only if $\ga\neq 0$ and its isometric sphere is then
$S_g=\{z\in\HH\;:\;\n(\gamma z+\delta)=1\}$ and $S^+_g=
\{(z,r)\in\hcr\;:\; \n(\gamma z+\delta)+r^2\geq 1\}$.

Recall that since $\Ga$ has finite covolume, every parabolic fixed
point $\xi$ of $\Ga$ is {\it bounded}, that is the quotient of
$\partial_\infty\HH^n_\RR-\{\xi\}$ by the stabilizer of $\xi$ in $\Ga$
is compact. Let $\D_\infty$ be a compact fundamental domain for the
action of the stabilizer of $\infty$ in $\Ga$ on $\RR^{n-1}$. Then
the {\it Ford fundamental domain} $\F_\Ga$ of $\Ga$ associated to
$\D_\infty$ is
$$
\F_\Ga=\Big(\bigcap_{g\in\Ga,\;g\infty\neq \infty}\;S^+_g\Big)
\cap\big(\D_\infty\times \;]0,+\infty[\big)\;.
$$
It is well known (see for instance \cite[page 239]{Beardon83}) that
$\F_\Ga$ is a fundamental domain for $\Ga$ acting on $\HH^n_\RR$ (in
particular, $\F_\Ga$ satisfies condition (i) of a weak fundamental
domain) and that the set of points at infinity of $\bigcap_{g\in\Ga,
  \;g\infty\neq \infty}\;S^+_g$ is a locally finite set of parabolic
fixed points in $\partial_\infty\HH^n_\RR$.  Furthermore, since
parabolic fixed points are bounded and have a precisely invariant
horoball centered at them, and since the tangency of a circle and its
tangent is quadratic, the condition (iii) is satisfied for every
$\epsilon$ small enough. Note that $\F_\Ga$ satisfies condition (ii)
with $K=\D_\infty$.

\bigskip Let us fix a weak fundamental domain $\F$ for the action of
$\SLO$ on $\hcr$. A positive definite form $f\in\Q^+(\OOO)$ is {\it
  reduced} if $\Phi([f])\in\F$ and an indefinite form
$f\in\Q^\pm(\OOO)$ is {\it reduced} if $\Psi([f])\cap\F\neq 0$. We say
that a negative definite form $f\in-\Q^+(\OOO)$ is {\it reduced} if
$-f$ is reduced. The notion of being reduced does depend on the choice
of a weak fundamental domain, which allows to choose it appropriately
when computing examples. Recall that $\Q(\Delta)$ is equal to
$\Q^\pm(\Delta)$ if $\Delta>0$ and to $\Q^+(\Delta)\cup-\Q^+(\Delta)$
if $\Delta<0$.

\btheo \label{theo:reducgeom}
For every $\Delta\in\ZZ-\{0\}$, the number of reduced elements of
$\Q(\OOO,\Delta)$ is finite. 
\etheo

Theorem \ref{theo:inducintro} in the Introduction is a restatement of
this one.

\medskip \dem Note that the Euclidean norm on $\HH$ is
$\|z\|=\n(z)^{\frac{1}{2}}$.

\medskip Let us first prove that the number of reduced elements of
$\Q^+(\OOO,\Delta)$ is finite.  

\medskip For every $f\in\Q^+(\OOO,\Delta)$, let $a=a(f)>0$, $b=b(f)$
and $c=c(f)$. We have $\n(b)-ac=\Delta<0$, hence $c$ is determined by
$a$ and $b$.  The form $f$ is reduced if and only if $\Phi([f])=
(-\frac{b}{a}, \frac{\sqrt{-\Delta}}{a}) \in \F$.  By the condition
(ii) and since $K$ is compact, $\|\frac{b}{a}\|$ is bounded.  Hence,
if we have an upper bound on $a$, by the discreteness of $\OOO$, the
elements $a$ and $b$ may take only finitely many values, and so does
$c$, therefore the result follows.

Let $\kappa,\epsilon,Z$ be as in the condition (iii).  If
$\frac{\sqrt{-\Delta}}{a}\geq \epsilon$, then $a$ is bounded from
above, and we are done. Otherwise, by condition (iii), there exists
$s$ in the finite set $Z$ such that $\Phi([f])\in\E_s$. In particular,
$$
\Big\|-\frac{b}{a}-s\,\Big\|\leq 
\kappa \,\Big(\frac{\sqrt{-\Delta}}{a}\Big)^2\;.
$$
Since the set of parabolic elements of $\SLO$ is $A\cup\{\infty\}$, we
may write $s=\frac{u}{v}$ with $u\in\OOO$ and $v\in\NN-\{0\}$. The
above inequality becomes
$$
a\;\|\,bv+au\,\|\leq \kappa\,|\Delta|\,v\;.
$$
The element $bv+au\in\OOO$ either is equal to $0$ or has reduced norm,
hence Euclidean norm, at least $1$. In the second case, we have an
upper bound on $a$, as wanted. In the first case, we have $\frac{b}{a}
= -\frac{u}{v}$, that is $b=-\frac{au}{v}$. Hence
$$ 
\Delta v^2=(\n(b)-ac)v^2=a(a\n(u)-cv^2)\;.
$$
Since $a\n(u)-cv^2\in\ZZ$, the integer $a$ divides the nonzero
integer $\Delta v^2$, hence $a$ is bounded, as wanted.

\medskip Let us now prove that the number of reduced elements of
$\Q^\pm(\OOO,\Delta)$ is finite, which concludes the proof of Theorem
\ref{theo:reducgeom}.

\medskip We have $\Delta >0$. With $K$ a compact subset as in the
condition (ii), let $\delta=\sup_{x\in K}\|x\|$. Let
$f\in\Q^\pm(\OOO,\Delta)$ be reduced, and fix $(z,r)\in \C(f)\cap
\F$. Let $a=a(f)$, $b=b(f)$ and $c=c(f)$.

Assume first that $a=0$. Then $\n(b)=\Delta$, hence $b$ takes only
finitely many values, by the discreteness of $\OOO$. Recalling that
$\C(f)=\{(z,r)\in \hcr\;:\;\tr(\overline{z}b)+c=0\}$, we have by
the Cauchy-Schwarz inequality
$$
|c|=|\tr(\overline{z}b)|\leq 2\,\|z\|\,\|b\|\leq 
2\,\delta\sqrt{\Delta}\;.
$$
Again by discreteness, $c$ takes only finitely many values, and the
result follows.

Assume that $a\neq 0$, and up to replacing $f$ by $-f$ (which is
reduced if $f$ is), that $a>0$.  We have $\n(b)-ac=\Delta$, hence $c$
is determined by $a$ and $b$. Recalling that $\C(f)=\{(z,r)\in
\hcr\;:\;\n(az+b)+a^2r^2= \Delta\}$, we have by the triangular
inequality
$$
\Big\|\frac{b}{a}\Big\|\leq \Big\|z+\frac{b}{a}\Big\|+\|z\|\leq
\sqrt{\Delta}+\delta\;.
$$ 
Hence as in the positive definite case, if we have an upper bound on
$a$, the result follows.

Let $\kappa,\epsilon,Z$ be as in the condition (iii). Note that $r\leq
\frac{\sqrt{\Delta}}{a}$. Hence if $r\geq \epsilon$, then we have an
upper bound $a\leq \frac{\sqrt{\Delta}}{\epsilon}$, as wanted.
Therefore, we may assume that $(z,r)$ belongs to $\C(f)\cap\E_s$ for
some $s\in Z$. In particular,
$$
\Big\|\,z+\frac{b}{a}\,\Big\|= \sqrt{\frac{\Delta}{a^2}-r^2}\;\;\;
{\rm and}\;\;\;
\|z-s\|\leq \kappa\,r^2\;.
$$

First assume that $\|\frac{b}{a}+s\|\geq \frac{\sqrt{\Delta}}{a}$.
Then by the inverse triangular inequality
$$
\kappa \,r^2\geq \|s-z\|\geq 
\Big\|\frac{b}{a}+s\Big\|-\Big\|z+\frac{b}{a}\,\Big\| 
\geq \frac{\sqrt{\Delta}}{a} -  \sqrt{\frac{\Delta}{a^2}-r^2}
\geq \frac{r^2}{2\frac{\sqrt{\Delta}}{a}}\;.
$$
Therefore, we have an upper bound $a\leq 2\,\kappa \sqrt{\Delta}$, as
wanted.

Now assume that $\|\frac{b}{a}+s\|< \frac{\sqrt{\Delta}}{a}$.  Write
$s=\frac{u}{v}$ with $u\in\OOO$ and $v\in\NN-\{0\}$. We have $\n(au
+bv) < \Delta v^2$. The element $w=au+bv$, belonging to $\OOO$ and
having reduced norm at most $\Delta v^2$, can take only finitely many
values. The  positive integer $v^2\Delta -\n(w)$
is equal to
$$
v^2(\n(b) -ac) -\n(au+bv)=-\tr(\overline{au}\,bv)-\n(au)-v^2ac
=-a(\tr(\overline{u}\,b\,v)+a\n(u)+v^2c)\;.
$$
Since $\tr(\overline{u}\,b\,v)+a\n(u)+v^2c\in\ZZ$ by the properties of
the reduced norm, the reduced trace and the conjugate of elements of
$\OOO$, this implies that the integer $a$ divides the nonzero integer
$v^2\Delta -\n(w)$, hence $a$ is bounded, as wanted.

This concludes the proof of Theorem \ref{theo:reducgeom}. \cqfd

\bcoro \label{coro:reduc}
For every $\Delta\in\ZZ-\{0\}$, the number of orbits of $\SLO$
in $\Q(\OOO,\Delta)$, hence in $\Q^+(\OOO,\Delta)$ and in
$\Q^\pm(\OOO,\Delta)$, is finite.
\ecoro

\dem This immediately follows from Theorem \ref{theo:reducgeom}, by
the equivariance properties in Proposition \ref{prop:equivar} and the
assumption (i) on a weak fundamental domain (that was not used in the
proof of Theorem \ref{theo:reducgeom}). \cqfd

\medskip\noindent {\bf Example. }  Let $A$ be Hamilton's quaternion
algebra over $\QQ$. Let $\OOO$ be Hurwitz's maximal order in $A$, and
let $\OOO'=\ZZ+\ZZ i+\ZZ j+\ZZ k$ be the order of Lipschitz integral
quaternions.

We identify $\HH$ and $\RR^4$ by the $\RR$-linear map sending
$(1,i,j,k)$ to the canonical basis of $\RR^4$.  Let $V\subset\OOO'$
denote the set of vertices of the $4$-dimensional unit cube $[0,1]^4$.
We claim that the set
$$
\F=\{(z,r)\in\HH^5_\RR:z\in[0,1]^4, \n(z-s)+r^2\ge 1 \text{ for all }
  s\in V\}
$$
is a weak fundamental domain for $\SL_2(\OOO')$, hence for $\SLO$.
For every $s\in V$, the $3$-sphere in $\HH$ with equation $\n(z-s)= 1$
is the isometric sphere of $\begin{pmatrix} 0& -1\\1 & \ \
  s \end{pmatrix}\in\SL_2(\OOO')$.  Since the diameter of the cube
$[0,1]^4$ is $2$, the closed balls bounded by these spheres cover
$[0,1]^4$. This unit cube is a fundamental polytope of the subgroup of
unipotent elements of $\SL_2(\OOO')$ fixing $\infty$.  Thus, $\F$
contains a Ford fundamental domain of $\SL_2(\OOO')$, which implies
property (i) of a weak fundamental domain. Property (ii) (with $K$ the
unit cube) is valid by the definition of $\F$. Property (iii) follows
from the fact that the only point at infinity of $\F$ besides $\infty$
is the center point $\frac{1+i+j+k}{2}$ of the unit cube, which is the
only point of this cube which does not belong to one of the open balls
whose boundary is one of the isometric spheres used to define
$\F$. Note that $\frac{1+i+j+k}{2}\in A$ is a parabolic fixed point of
$\SL_2(\OOO')$.

Recall that a positive definite Hamiltonian form
$f\in\Q^+(\OOO,\Delta)$ with coefficients $a=a(f)$,
$b=b(f)=b_1+b_2i+b_3j+b_4k$ and $c=c(f)$ is reduced (for this choice
of weak fundamental domain) if $(-\frac
ba,\frac{\sqrt{-\Delta}}a)\in\F$. A straightforward manipulation of
the defining inequalities of $\F$ shows that $f\in\Q^+(\OOO,\Delta)$
is reduced if and only if its coefficients satisfy the following set
of $25$ inequalities
\begin{equation}\label{eq:reduced}
  a   >0\,,\quad \quad
  0  \le-b_\ell\le a\,, \quad \quad
  a\big(a-c -2\sum_{m\in P} b_m\big)\le\card(P)
\end{equation}
for all $\ell\in\{1,2,3,4\}$ and for all subsets
$P\subset\{1,2,3,4\}$.  Theorem \ref{theo:reducgeom} implies that
there are only a finite number of forms in $\Q^+(\OOO,\Delta)$ whose
coefficients satisfy the inequalities \eqref{eq:reduced}.

\medskip Similarly, an indefinite Hamiltonian form
$f\in\Q^\pm(\OOO,\Delta)$ with $a(f)=a>0$, $b(f)=b_1+b_2i+b_3j+b_4k$
and $c(f)=c$ is reduced, that is $\C(f)$ meets $\F$, if and only if
the following system of $16$ linear inequalities and one quadratic
inequality in four real variables $X_1,X_2,X_3,X_4$ has a solution in
the unit cube $[0,1]^4$:
$$
\sum_{\ell=1}^4\; 2X_\ell\frac{b_\ell}{a} +{X_\ell}^2\leq -\frac{c}{a}
\;,\quad\quad\quad
\sum_{\ell=1}^4\; 2X_\ell\frac{b_\ell}{a}+\sum_{m\in P}\;2X_m
\leq -1 -\frac{c}{a} +\card(P)\;,
$$
for all subsets $P\subset\{1,2,3,4\}$.

\appendix
\section{The hyperbolic covolume of $\SLO$, by Vincent Emery}
\label{sect:hypercovol}

Let $A$ be a definite quaternion algebra over $\QQ$, with reduced
discriminant $D_A$, and let $\OOO$ be a maximal order in $A$ (see for
instance \cite{Vigneras80} and Section \ref{sect:defquadalg} for
definitions and properties). Given a quaternion algebra $A'$ over a
field $k$, let $\SL_2(A')=\SL_1(\M_2(A'))$ be the group of elements of
the central simple $2\times 2$ matrix algebra $\M_2(A')$ having reduced
norm $1$. For any subring $\OOO'$ of $A'$, let
$\SL_2(\OOO')=\SL_2(A')\cap\M_2(\OOO')$ and
$\PSL_2(\OOO')=\SL_2(\OOO')/\{\pm \id\}$. Fixing an identification
between $A\otimes_\QQ \RR$ and Hamilton's real quaternion algebra
$\HH$ turns $\SLO$ into an arithmetic lattice in $\SL_2(\HH)$. Hence
$\SLO$ acts by isometries with finite covolume on the real hyperbolic
space $\hcr$ (see for instance Section \ref{sect:hambiagroup} for
generalities).

In this appendix, the following result is proved using Prasad's volume
formula in \cite{Prasad89}. See the main body of this paper for a
proof using Eisenstein series.  The author of this appendix thanks
Jouni Parkkonen and Fr\'ed\'eric Paulin for helpful discussions.  He
is particularly grateful to Fr\'ed\'eric Paulin for his help concerning the
Bruhat-Tits buildings appearing in the proof. This work is supported
by the Swiss National Science Foundation, Project number {\tt
  PP00P2-128309/1}.

\btheo\label{theo:appendix}
The hyperbolic covolume of $\SLO$ is
$$
\covol(\SLO)=\frac{\zeta(3)}{11520}\prod_{p\mid D_A} (p^3-1)(p-1)\;,
$$
where $p$ ranges over the prime integers.
\etheo

\dem Let $\P$ be the set of positive primes in $\ZZ$. For every $p\in
\P$, let $\OOO_p=\OOO\otimes_\ZZ\ZZ_p$, which is a maximal order in
the quaternion algebra $A_p=A\otimes_\QQ \QQ_p$ over $\QQ_p$ (see for
instance \cite[page 84]{Vigneras80}).

We refer for instance to \cite{Tits66} for the classification of the
semi-simple connected algebraic groups over $\QQ$. Let ${\bf G}$ be
the (affine) algebraic group over $\QQ$, having as group of
$K$-points, for each characteristic zero field $K$, the group
$$
{\bf G}(K)=\SL_2(A\otimes_\QQ K)=\SL_1(\M_2(A\otimes_\QQ K))\;.
$$
The group ${\bf G}$ is absolutely (quasi-)simple and simply connected.
Indeed,  the $\CC$-algebra $A\otimes_\QQ \CC$ is isomorphic to
$\M_2(\CC)$ and thus the complex Lie group ${\bf G}(\CC)$ is isomorphic
to $\SL_1(\M_4(\CC))=\SL_4(\CC)$ (note that we are using the reduced
norm and not the norm). Furthermore, ${\bf G}$ is an inner form of the
split algebraic group $\G=\SL_4$ over $\QQ$. The (absolute)
rank of $\G$ and the exponents of $\G$ are given by
\begin{equation}\label{eq:rankexpo}
r=3 \;\;\;{\rm and}\;\;\;m_1=1,\; m_2=2,\;m_3=3
\end{equation}
(see for instance \cite[page 96]{Prasad89}). We consider the
$\ZZ$-form of ${\bf G}$ such 
that ${\bf G}(\ZZ)=\SLO$ and ${\bf G}
(\ZZ_p)=\SL_2(\OOO_p)$ for every $p\in\P$ (see for instance \cite[page
382]{ParPau10GT} for details).

Let $\I_{{\bf G},\QQ_p}$ be the Bruhat-Tits building of ${\bf G}$ over
$\QQ_p$ (see for instance \cite{Tits79} for the necessary background
on Bruhat-Tits theory).  Recall that a subgroup of ${\bf G}(\QQ_p)$ is
{\it parahoric} if it is the stabilizer of a simplex of $\I_{{\bf
    G},\QQ_p}$; a {\it coherent family of parahoric subgroups} of
${\bf G}$ is a family $(Y_p)_{p\in\P}$, where $Y_p$ is a parahoric
subgroup of ${\bf G} (\QQ_p)$ and $Y_p= {\bf G} (\ZZ_p)$ for $p$ big
enough. The {\it principal lattice} associated with this family is the
subgroup ${\bf G}(\QQ)\cap \prod_{p} Y_p$ of ${\bf G}(\QQ)$
(diagonally contained in the group ${\bf G}(\AA_f)=\prod_{p}'{\bf
  G}(\QQ_p)$ of finite ad\`eles of ${\bf G}$, where as usual $\prod'$
indicates the restricted product).

For every $p\in\P$, recall that by the definition of the discriminant
$D_A$ of $A$, if $p$ does not divide $D_A$, then the algebra $A_p$ is
isomorphic to $\M_2(\QQ_p)$, and otherwise $A_p$ is a
$d^2$-dimensional central division algebra with center $\QQ_p$ with
$d=2$. Furthermore, for the discrete valuation $\nu=\nu_p\circ\n$
where $\nu_p$ is the discrete valuation of $\QQ_p$ and $\n$ the
reduced norm on $A_p$, the maximal order $\OOO_p$ is equal to the
valuation ring of $\nu$ (see for instance \cite[page
34]{Vigneras80}). 

First assume that $p$ does not divide $D_A$. Then ${\bf G}$ is
isomorphic to $\G=\SL_4$ over $\QQ_p$. The vertices of the building
$\I_{{\bf G} ,\QQ_p}$ are the homothety classes of $\ZZ_p$-lattices in
${\QQ_p}^4$. In particular $\SL_2(\OOO_p)=\SL_4(\ZZ_p)$ is the
stabilizer of the class of the standard $\ZZ_p$-lattice ${\ZZ_p}^4$,
hence is parahoric.

Now assume that $p$ divides $D_A$. Then ${\bf G}(\QQ_p)= \SL_m(A_p)$
with $m=2$ and ${\bf G}(\QQ_p)$ has local type $^dA_{md-1}=\;^2A_3$ in
Tits' classification \cite[\S 4.4]{Tits79}.  The corresponding local
index is shown below:
\begin{center}
\includegraphics[width=3cm]{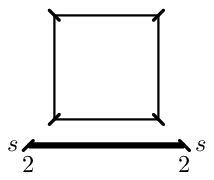}
\end{center}
\begin{center} Local index of type $^2A_3$
\end{center}
The building $\I_{{\bf G}, \QQ_p}$ is a tree (see for instance
\cite{Serre83} for the construction of the Bruhat-Tits tree of
$\SL_2(K)$ even when $K$ is a non-commutative division algebra endowed
with a discrete valuation). Its vertices are the homothety classes of
$\OOO_p$-lattices in the right $A_p$-vector space ${A_p}^2$. In
particular $\SL_2(\OOO_p)$ is the stabilizer of the class of the
standard $\OOO_p$-lattice ${\OOO_p}^2$, hence is parahoric.

Therefore, by definition the family $(\SL_2(\OOO_p))_{p\in\P}$ is a
coherent family of (maximal) parahoric subgroups of ${\bf G}$, and
$\SLO={\bf G}(\ZZ)={\bf G}(\QQ)\cap \prod_{p\in\P} {\bf G} (\ZZ_p)$ is
its associated principal lattice.

\medskip For every $p\in\P$, let $\ov M_p$ (respectively $\ov\M_p$) be
the Levi subgroup, defined over the residual field $\FF_p=
\ZZ_p/p\ZZ_p$, of the identity component of the reduction modulo $p$
of the smooth affine group scheme over $\ZZ_p$ associated with the
vertex of $\I_{{\bf G}, \QQ_p}$ (respectively $\I_{\G,\QQ_p}$)
stabilized by the parahoric subgroup $\SL_2(\OOO_p)$ (respectively
$\SL_4(\ZZ_p)$); see for instance \cite[\S 3.5]{Tits79}.  Note that
$\ov M_p=\ov\M_p$ if $p$ does not divide $D_A$, and that for every
$p\in\P$ the algebraic group $\ov\M_p$ is isomorphic to $\SL_4$ over
$\FF_p$. In particular $\ov\M_p(\FF_p)=\SL_4(\FF_p)$ and  thus, for
every $p\in\P$,
\begin{equation}\label{eq:dimsplit}
\dim \ov\M_p=15 \;\;{\rm and}\;\;
|\ov \M_p(\FF_p)|=p^6(p^2-1)(p^3-1)(p^4-1)\;,
\end{equation}
the orders of finite groups of Lie type being listed for example
in~\cite[Table~1]{Ono66}.  If $p$ divides $D_A$, by applying the theory
in~\cite[\S 3.5.2]{Tits79} on the local index $^2A_3$, we see that the
semi-simple part $\ov M_p^{ss}$ of $\ov M_p$ (given as the commutator
algebraic group $[ \ov M_p, \ov M_p]$) is of type $^2(A_1\times A_1)$
and the radical $R(\ov M_p)$ of $\ov M_p$ must be a one-dimensional
non-split torus over $\FF_p$. In particular $|R(\ov M_p)(\FF_p)| = p +
1$ and $\ov M_p^{ss}(\FF_p)$ has the same order as $\SL_2(\FF_{p^2})$,
that is $p^2 (p^4 - 1)$. Since the radical $R(\ov M_p)$ is
central in $\ov M_p$ and the intersection $R(\ov M_p) \cap \ov
M_p^{ss}$ is finite (see~\cite[Prop.~7.3.1]{Springer98}), the product map
$$
	\begin{array}[h]{ccc}
	\ov M_p^{ss} \times R(\ov M_p) &\to& \ov M_p\\
	(x,y) &\mapsto& xy	
	\end{array}
	\label{eq:isogeny-2A3}
$$
is an isogeny (defined over $\FF_p$) and using Lang's isogeny theorem
(see for example~\cite[Prop~6.3, page~290]{PlaRap94}), we obtain the
order of $\ov M_p(\FF_p)$ as the product $|\ov M_p^{ss}(\FF_p)| \cdot
|R(\ov M_p)(\FF_p)|$. Therefore, for every $p\in\P$ dividing $D_A$,
\begin{equation}\label{eq:dimnonsplit}
\dim \ov M_p=7 \;\;{\rm and}\;\;
|\ov M_p(\FF_p)|=p^2(p^4-1)(p+1)\;.
\end{equation}

\medskip Let $\mu$ be the Haar measure on ${\bf G}(\RR)=\SLH$
normalized as in \cite{Prasad89}. That is, if $w$ is the top degree
exterior form on the real Lie algebra of ${\bf G}(\RR)$ whose
associated invariant differential form on ${\bf G}(\RR)$ defines the
measure $\mu$, if ${\bf G}_u(\RR)$ is a compact real form of ${\bf
  G}(\CC)$, then the complexification $w_\CC$ of $w$ on the complex
Lie algebra of ${\bf G} (\CC)={\bf G}_u(\CC)$ defines a top degree
exterior form $w_u$ on the real Lie algebra of ${\bf G}_u(\RR)$, whose
associated invariant differential form on ${\bf G}_u(\RR)$ defines a
measure $\mu_u$, and we require that $\mu_u({\bf G}_u(\RR))=1$.

Let $\mu'$ be the Haar measure on $\PSL_2(\HH)=
\operatorname{SO}_0(1,5)$ which disintegrates by the fibration
$\operatorname{SO}_0(1,5)\ra \operatorname{SO}_0(1,5) /
\operatorname{SO}(5) =\HH^5_\RR$, with measures on the fibers of total
mass one $1$, and measure on the base the Riemannian measure
$\dvol_{\HH^5_\RR}$ of the Riemannian metric of constant sectional
curvature $-1$.  Let $\wt\mu'$ be the Haar measure on $\SL_2(\HH)$
such that the tangent map at the identity of the double cover of real
Lie groups $\SL_2(\HH)\ra \PSL_2(\RR)$ preserves the top degree
exterior forms defining the Haar measures. In particular, since $-\id$
belongs to $\SLO$,
\begin{align}
\covol(\SLO)&=\Vol_{\HH^5_\RR}(\PSLO\bs\HH^5_\RR)= 
\mu'(\PSLO\bs\PSL_2(\HH))\nonumber\\ 
\label{eq:decovolahaar}&=\wt\mu'(\SLO\bs\SLH)\;.
\end{align}
Similarly, with $\SSS_5$ the $5$-sphere endowed with its standard
Riemannian metric of constant sectional curvature $+1$, let $\mu'_u$
be the Haar measure on $\operatorname{SO}(6)$ which disintegrates by
the fibration $\operatorname{SO}(6)\ra \operatorname{SO}(6) /
\operatorname{SO}(5)= \SSS_5$ with measures on the fibers of total
mass one $1$, and measure on the base the Riemannian measure. In
particular, $\mu'_u(\operatorname{SO}(6))=\Vol(\SSS_5)$. Recall that
$\Vol(\SSS_n)= \frac{2\pi^{m}}{(m-1)!}$ if $n=2m-1\geq 3$. It is well
known (see for instance \cite{Helgason78}) that the duality
$G/K\mapsto G_u/K$ between irreducible symmetric spaces of non-compact
type endowed with a left invariant Riemannian metric and the ones of
compact type, where $G_u$ is a compact form of the complexification of
$G$, sends $\HH^5_\RR$ on $\SSS_5$, hence $\mu'$ on $\mu'_u$.

The maximal compact subgroup $\operatorname{SU}(4)$ of $\SL_4(\CC)$
is a covering of degree $2$ of $\operatorname{SO}(6)$, which is the
compact real form corresponding to $\operatorname{SO}_0(1,5)$.  Hence we
have (as first proved in \cite[\S 13.3]{Emery09})
\begin{equation}\label{eq:comparHaar}  
\wt\mu'=2\Vol(\SSS_5)\;\mu= 2\pi^3\;\mu\;.
\end{equation}

By Prasad's volume formula \cite[Theo.~3.7]{Prasad89} (where with the
notation of this theorem, $\ell=k=\QQ$ hence $D_k=D_\ell=1$,
$S=V_\infty= \{\infty\}$ and the Tamagawa number $\tau_\QQ(\bf G)$ is
$1$), we have, since $\ov\M_p=\ov M_p$ if $p$ does not divide $D_A$
and by Equation \eqref{eq:rankexpo} for the second equality,
\begin{align}
\mu(\SLO\bs\SLH)&=\prod_{i=1}^r
\frac{(m_i)!}{(2\pi)^{m_i+1}}\prod_{p\in\P}
\frac{p^{(\dim \ov M_p+\dim
    \ov \M_p)/2}}{|\ov M_p(\FF_p)|}\nonumber\\ 
\label{eq:formulprasad}&=
\frac{12}{(2\pi)^9}\prod_{p\in\P}\frac{p^{\dim \ov \M_p}}{|\ov
    \M_p(\FF_p)|}\;
\prod_{p|D_A}\frac{|\ov \M_p(\FF_p)|}{|\ov M_p(\FF_p)|}\;
p^{(\dim \ov M_p-\dim \ov \M_p)/2}.
\end{align}
Using Euler's product formula
$\zeta(s)=\prod_{p\in\P}\frac{1}{1-p^{-s}}$ for Riemann's zeta
function, we have by Equation \eqref{eq:dimsplit}, since
$\zeta(2)=\frac{\pi^2}{6}$ and $\zeta(4)=\frac{\pi^4}{90}$,
\begin{equation}\label{eq:prodzeta}  
\prod_{p\in\P}\frac{p^{\dim \ov \M_p}}{|\ov
    \M_p(\FF_p)|}=\zeta(2)\zeta(3)\zeta(4)=\frac{\pi^6\zeta(3)}{540}\;.
\end{equation}
Using the equations
\eqref{eq:decovolahaar}, \eqref{eq:comparHaar}, \eqref{eq:formulprasad},
\eqref{eq:prodzeta}, \eqref{eq:dimsplit} and \eqref{eq:dimnonsplit},
the result follows. \cqfd

{\small \bibliography{../biblio} }

\begin{thebibliography}{EGM}

\bibitem[Asl]{Aslaksen96}
H.~Aslaksen.
\newblock {\it Quaternionic determinants}.
\newblock {Math.~Intelligencer {\bf 18} (1996) 57--65}.

\bibitem[Bea]{Beardon83}
A.~F. Beardon.
\newblock {\it The geometry of discrete groups}.
\newblock {Grad. Texts Math. {\bf 91}, Springer-Verlag, 1983}.

\bibitem[BHC]{BorelHarishChandra62}
A.~Borel and Harish-Chandra.
\newblock {\it Arithmetic subgroups of algebraic groups}.
\newblock {Ann. of Math. {\bf 75} (1962) 485--535}.

\bibitem[Bor]{Borel66}
A.~Borel.
\newblock {\it Reduction theory for arithmetic groups}.
\newblock {in ``Algebraic Groups and Discontinuous Subgroups'', A.~Borel and
  G.~D.~Mostow eds, Proc. Sympos. Pure Math. (Boulder, 1965), pp. 20--25, Amer.
  Math. Soc. 1966}.

\bibitem[Bou]{Bourbaki59}
N.~Bourbaki.
\newblock {\it Algèbre. Chapitre 9 : Formes sesquilinéaires et formes
  quadratiques}.
\newblock {Hermann, 1959}.

\bibitem[BH]{BreHelm96}
S.~Breulmann and V.~Helmke.
\newblock {\it The covolume of quaternion groups on the four-dimensional
  hyperbolic space}.
\newblock {Acta Arith. {\bf 75} (1996) 9--21}.

\bibitem[BT]{BruTit72}
F.~Bruhat and J.~Tits.
\newblock {\it Groupes r\'eductifs sur un corps local, I: Donn\'ees radicielles
  valu\'ees}.
\newblock {Pub. Math. I.H.\'E.S. {\bf 41} (1972) 5--251}.

\bibitem[BV]{BucVol07}
J.~Buchmann and U.~Vollmer.
\newblock {\it Binary quadratic forms : an algorithmic approach}.
\newblock {Alg.~Comput.~Math. {\bf 20}, Springer Verlag, 2007}.

\bibitem[Cas]{Cassels08}
J.~W.~S. Cassels.
\newblock {\it Rational quadratic forms}.
\newblock {Dover, 2008}.

\bibitem[Deu]{Deuring35}
M.~Deuring.
\newblock {\it Algebren}.
\newblock {2nd ed., Erg. Math. Grenz. {\bf 41}, Springer Verlag, 1968}.

\bibitem[Die]{Dieudonne71}
J.~Dieudonn\'e.
\newblock {\it Les déterminants sur un corps non commutatif}.
\newblock {Bull. Soc. Math. France, {\bf 71} (1943) 27--45}.

\bibitem[Eic]{Eichler38}
M.~Eichler.
\newblock {\it Über die Idealklassenzahl total definiter
  Quaternionenalgebren}.
\newblock {Math. Z. {\bf 43} (1938) 102--109}.

\bibitem[EGM]{ElsGruMen98}
J.~Elstrodt, F.~Grunewald, and J.~Mennicke.
\newblock {\it Groups acting on hyperbolic space: Harmonic analysis and number
  theory}.
\newblock {Springer Mono. Math., Springer Verlag, 1998}.

\bibitem[Eme]{Emery09}
V.~Emery.
\newblock {\it Du volume des quotients arithmétiques de l'espace
  hyperbolique}.
\newblock {Thèse n$^0$ 1648, Univ. Fribourg (Suisse), 2009}.

\bibitem[ERS]{EskRudSar91}
A.~Eskin, Z.~Rudnick, and P.~Sarnak.
\newblock {\it A proof of Siegel's weight formula}.
\newblock {Internat. Math. Res. Notices {\bf 5} (1991) 65--69}.

\bibitem[HI]{HasIbu80}
K.~Hashimoto and T.~Ibukiyama.
\newblock {\it On class numbers of positive definite binary quaternion
  Hermitian forms. I, II, III}.
\newblock {J. Fac. Sci. Univ. Tokyo {\bf 27} (1980) 549--601, {\bf 28} (1981)
  695--699, {\bf 30} (1983) 393--401}.

\bibitem[Hel]{Helgason78}
S.~Helgason.
\newblock {\it Differential geometry, Lie groups and symmetric spaces}.
\newblock {Academic Press, 1978}.

\bibitem[Hil]{Hild07}
T.~Hild.
\newblock {\it The cusped hyperbolic orbifolds of minimal volume in dimensions
  less than ten}.
\newblock {J. Algebra {\bf 313} (2007) 208--222}.

\bibitem[JW]{JohWei99}
N.~W. Johnson and A.~I. Weiss.
\newblock {\it Quaternionic modular groups}.
\newblock {Linear Algebra Appl. {\bf 295} (1999) 159--189}.

\bibitem[Kel]{Kellerhals03}
R.~Kellerhals.
\newblock {\it Quaternions and some global properties of hyperbolic
  $5$-manifolds}.
\newblock {Canad. J. Math. {\bf 55} (2003) 1080--1099}.

\bibitem[KO]{KraOse90}
V.~Krafft and D.~Osenberg.
\newblock {\it Eisensteinreihen f\"ur einige arithmetisch definierte
  Untergruppen von $\SL_2(\HH)$}.
\newblock {Math. Z. {\bf 204} (1990) 425--449}.

\bibitem[Lan]{Langlands66b}
R.~P. Langlands.
\newblock {\it The volume of the fundamental domain for some arithmetical
  subgroups of Chevalley groups}.
\newblock {in "Algebraic Groups and Discontinuous Subgroups" (Boulder, 1965)
  pp.~143–148, (Proc. Sympos. Pure Math. IX, Amer. Math. Soc. 1966}.

\bibitem[MOS]{MagObeSon66}
W.~Magnus, F.~Oberhettinger, and R.~Soni.
\newblock {\it Formulas and theorems for the special functions of mathematical
  physics}.
\newblock {3rd ed, Grund. math. Wiss. {\bf 52}, Springer Verlag, 1966}.

\bibitem[Ono]{Ono66}
T.~Ono.
\newblock {\it On algebraic groups and discontinuous groups}.
\newblock {Nagoya Math. J. {\bf 27} (1966) 279--322}.

\bibitem[PP1]{ParPau10GT}
J.~Parkkonen and F.~Paulin.
\newblock {\it Prescribing the behaviour of geodesics in negative curvature}.
\newblock {Geometry \& Topology {\bf 14} (2010) 277-392}.

\bibitem[PP2]{ParPau}
J.~Parkkonen and F.~Paulin.
\newblock {\it \'Equidistribution, comptage et approximation par irrationnels
  quadratiques}.
\newblock { J. Mod. Dyn.   {\bf  6} (2012) 1-40}.

\bibitem[PP3]{ParPau11CJM}
J.~Parkkonen and F.~Paulin.
\newblock {\it On the representations of integers by indefinite binary
  Hermitian forms}.
\newblock {Bull. London Math. Soc.  {\bf 43} (2011), 1048-1058}.

\bibitem[PR]{PlaRap94}
V.~Platonov and A.~Rapinchuck.
\newblock {\it Algebraic groups and number theory}.
\newblock {Academic Press, 1994}.

\bibitem[Pra]{Prasad89}
G.~Prasad.
\newblock {\it Volumes of $S$-arithmetic quotients of semi-simple groups}.
\newblock {Publ. Math. Inst. Hautes \'Etudes Sci. {\bf 69} (1989) 91--117}.

\bibitem[Pro]{Pronin67}
L.~Pronin.
\newblock {\it Integral binary Hermitian forms over the skewfield of
  quaternions}.
\newblock {(Russian) Vestnik Har'kov. Gos. Univ. {\bf 26} (1967) 27--41}.

\bibitem[Ran]{Rankin39}
R.~A. Rankin.
\newblock {\it Contributions to the theory of Ramanujan’s function $\tau(n)$
  and similar arithmetical functions. I. The zeros of the function
  $\sum_{n=1}^\infty \tau(n)/n^s$ on the line $\Re s = 13/2$. II. The order of
  the Fourier coefficients of integral modular forms}.
\newblock {Proc. Cambridge Philos. Soc. {\bf 35} (1939) 351--372}.

\bibitem[Rei]{Reiner75}
I.~Reiner.
\newblock {\it Maximal orders}.
\newblock {Academic Press, 1972}.

\bibitem[Sar]{Sarnak83}
P.~Sarnak.
\newblock {\it The arithmetic and geometry of some hyperbolic three-manifolds}.
\newblock {Acta Math. {\bf 151} (1983) 253--295}.

\bibitem[Sch]{Schoeneberg39}
B.~Schoeneberg.
\newblock {\it \"{U}ber die $\zeta$-Funktion einfacher hyperkomplexer Systeme}.
\newblock {Math. Ann. {\bf 117} (1939) 85--88}.

\bibitem[Sel]{Selberg40}
A.~Selberg.
\newblock {\it Bemerkungen \"uber eine Dirichletsche Reihe, die mit der Theorie
  der Modulformen nahe verbunden ist}.
\newblock {Arch. Math. Naturvid. {\bf 43} (1940) 47--50}.

\bibitem[Ser]{Serre83}
J.-P. Serre.
\newblock {\it Arbres, amalgames, SL$_2$}.
\newblock {3ème éd. corr., Ast\'erisque {\bf 46}, Soc. Math. France, 1983}.

\bibitem[Spr]{Springer98}
T.~A. Springer.
\newblock {\it Linear algebraic groups}.
\newblock {2nd ed., Progr. in Math. {\bf 9}, Birkh\"auser, 1998}.

\bibitem[Tit1]{Tits66}
J.~Tits.
\newblock {\it Classification of algebraic semisimple groups}.
\newblock {in "Algebraic Groups and Discontinuous Subgroups" (Boulder, 1965)
  pp. 33–62, Proc. Symp. Pure Math. IX, Amer. Math. Soc. 1966}.

\bibitem[Tit2]{Tits79}
J.~Tits.
\newblock {\it Reductive groups over local fields}.
\newblock {in "Automorphic forms, representations and $L$-functions"
  (Corvallis, 1977), Part 1, pp. 29–69, Proc. Symp. Pure Math. XXXIII, Amer.
  Math. Soc., 1979}.

\bibitem[Vig]{Vigneras80}
M.~F. Vign{\'e}ras.
\newblock {\it Arithm\'etique des alg\`ebres de quaternions}.
\newblock {Lect. Notes in Math. {\bf 800}, Springer Verlag, 1980}.

\bibitem[Wey]{Weyl40}
H.~Weyl.
\newblock {\it Theory of reduction for arithmetical equivalence I, II}.
\newblock {Trans.~Amer.~Math.~Soc. {\bf 48} (1940) 126--164, {\bf 51} (1942)
  203--231}.

\bibitem[Zag]{Zagier81}
D.~Zagier.
\newblock {\it Zetafunktionen und quadratische Körper}.
\newblock {Hochschultext, Springer Verlag, 1981}.

\end{thebibliography}

\bigskip
{\small\noindent \begin{tabular}{l} 
Department of Mathematics and Statistics, P.O. Box 35\\ 
40014 University of Jyv\"askyl\"a, FINLAND.\\
{\it e-mail: parkkone@maths.jyu.fi}
\end{tabular}
\medskip

\noindent \begin{tabular}{l}
D\'epartement de math\'ematique, UMR 8628 CNRS, B\^at.~425\\
Universit\'e Paris-Sud 11,
91405 ORSAY Cedex, FRANCE\\
{\it e-mail: frederic.paulin@math.u-psud.fr}
\end{tabular}
\medskip

\noindent \begin{tabular}{l}
Section de mathématiques, 2-4 rue du Lièvre,\\
Case postale 64, 1211 Genève 4, SWITZERLAND \\
{\it e-mail: Vincent.Emery@unige.ch}
\end{tabular}
}

\end{document}